\documentclass{amsart}
\usepackage{amsmath,amsxtra,amssymb}
\usepackage[all]{xy}

\newcommand{\A}{{\mathcal A}}
\newcommand{\B}{{\mathcal B}}

\newcommand{\U}{{\mathcal U}}

\newtheorem{theorem}{Theorem}[section]
\newtheorem{cor}[theorem]{Corollary}
\newtheorem{prop}[theorem]{Proposition}
\newtheorem{rem}[theorem]{Remark}
\newtheorem{lemma}[theorem]{Lemma}

\begin{document}
\title{On the best constants in noncommutative Khintchine-type inequalities}
\author{Uffe Haagerup$^{(1)}$ and Magdalena Musat$^{(2)}$}
\address{$^{(1)}$ Department of Mathematics and
Computer Science, University of Southern Denmark, Campusvej 55, 5230
Odense M, Denmark.\\
$^{(2)}$Department of Mathematical Sciences, University of Memphis,
373 Dunn Hall, Memphis, TN, 38152, USA.}
\email{$^{(1)}$haagerup@imada.sdu.dk\\$^{(2)}$mmusat@memphis.edu}.

%\author{Magdalena Musat}
%\address{Department of Mathematics, 0112\\
%University of California, San Diego\\
%La Jolla, CA 92093-0112}
%\email{mmusat@math.ucsd.edu}
\date{}


\begin{abstract}
We obtain new proofs with improved constants of the Khintchine-type
inequality with matrix coefficients in two cases. The first case is
the Pisier and Lust-Piquard noncommutative Khintchine inequality for
$p=1$\,, where we obtain the sharp lower bound of $\frac1{\sqrt{2}}$
in the complex Gaussian case and for the sequence of functions
$\{e^{i2^nt}\}_{n=1}^\infty$\,. The second case is Junge's recent
Khintchine-type inequality for subspaces of the operator space
$R\oplus C$\,, which he used to construct a cb-embedding of the
operator Hilbert space $OH$ into the predual of a hyperfinite
factor. Also in this case, we obtain a sharp lower bound of
$\frac1{\sqrt{2}}$\,. As a consequence, it follows that any subspace
of a quotient of $(R\oplus C)^*$ is cb-isomorphic to a subspace of
the predual of the hyperfinite factor of type III$_1$\,, with
cb-isomorphism constant $\leq \sqrt{2}$\,. In particular, the
operator Hilbert space $OH$ has this property.
\end{abstract}

\maketitle
\section{Introduction}

Let $r_n(t)=\text{sgn}(\sin(2^nt\pi))$\,, $n\in \mathbb{N}$ denote
the Rademacher functions on $[0, 1]$\,. The classical Khintchine
inequality states that for every $0<p< \infty$\,, there exist
constants $A_p$ and $B_p$ such that
\begin{equation}\label{eq200000209} A_p\left(\sum\limits_{k=1}^n a_k^2\right)^{\frac12}\leq
\left(\int\limits_0^1 \left|\sum\limits_{k=1}^n a_k r_k\right|
^p\,dt\right)^{\frac1p}\leq
B_p\left(\sum\limits_{k=1}^na_k^2\right)^{\frac12}\,,
\end{equation} for arbitrary $n\in \mathbb{N}$ and $a_1\,, \ldots
\,,a_n\in \mathbb{R}$\,.

Suppose $A_p$ and $B_p$ denote the best constants for which
(\ref{eq200000209}) holds. While it is elementary to prove that
$B_p=1$ for $0< p\leq 2$ and $A_p=1$ for $2\leq p< \infty$\,, it
took the work of many mathematicians to settle all the other cases,
including Szarek \cite{Sz} who proved that $A_1=\frac1{\sqrt{2}}$
(thus solving a long-standing conjecture of Littlewood), Young
\cite{Yo} who computed $B_p$ for $p\geq 3$\,, and the first-named
author (cf. \cite{Ha2}) who computed $A_p$ and $B_p$ in the
remaining cases.

The Khintchine inequality and its generalization to certain classes
of Banach spaces are deeply connected with the study of the geometry
of those Banach spaces (see \cite{MaP}). Noncommutative
generalizations of the classical Khintchine inequality to the case
of matrix-valued coefficients were first proved by Lust-Piquard
\cite {Lus} in the case $1< p< \infty$\,, and by Pisier and
Lust-Piquard \cite{LPP} for $p=1$\,. Their method of proof follows
the classical harmonic analysis approach of deriving Khintchine
inequality for the sequence $\{e^{i2^nt}\}_{n=1}^\infty$ from a
Paley inequality, for which they proved a noncommutative version
(see Theorem II.1 in \cite{LPP}). As a consequence, the following
noncommutative Khintchine inequality holds (see Corollary II.2 in
\cite{LPP}). Given $d\,, n\in \mathbb{N}$ and $x_1, \ldots, x_d\in
M_n(\mathbb{C})$\,, then
\begin{eqnarray}\label{eq2000002091}
{\frac1{1+\sqrt{2}}}|||\{x_j\}_{j=1}^d|||^*&\leq
&\left\|\sum\limits_{j=1}^d x_j\otimes e^{i2^nt}\right\|_{L^1([0,
1]; S_1^n)}\leq |||\{x_j\}_{j=1}^d|||^*\,,
\end{eqnarray}
where, by definition,
\begin{eqnarray}\label{eq55565676}
|||\{x_i\}_{i=1}^d|||^*:&= & \inf \left\{
\text{Tr}\left(\left(\sum\limits_{i=1}^d y_i^*y_i\right)^{\frac12} +
\left(\sum\limits_{i=1}^d z_iz_i^*\right)^{\frac12} \right);
x_i=y_i+z_i\in M_n(\mathbb{C})\right\}\,.
\end{eqnarray}
Here $S_1^n$ is $M_n(\mathbb{C})$ with the norm
$\|x\|_1:=\text{Tr}((x^*x)^{1/2})$\,, and $\text{Tr}$ is the
non-normalized trace on $M_n(\mathbb{C})$\,. We should also point
out that it was noted in the paper \cite{LPP} (cf. p. 250) that, by
using the lacunary sequence $\{3^n\}_{n\geq 1}$ instead of the
sequence $\{2^n\}_{n\geq 1}$, the lower bound in the inequality
(\ref{eq2000002091}) can be improved to $\frac12$\,.

By classical arguments (cf. Proposition 3.2 in \cite{Pi4})\,, if one
replaces $\{e^{i2^nt}\}_{n=1}^\infty$ by a sequence of independent
complex Gaussian, respectively, Rademacher or Steinhauss random
variables, the corresponding Khintchine inequality with matrix
coefficients follows, as well, with possibly different constants.

Our method, leading to improved constants, was inspired by ideas of
Pisier from \cite{Pi5}, and it is based on proving first directly
the dual inequality to (\ref{eq2000002091}) with constant
$\sqrt{2}$\,, where $\{e^{i2^nt}\}_{n=1}^\infty$ is replaced by a
sequence of independent complex-valued standard Gaussian random
variables on some probability space $(\Omega, \mathbb{P})$\,. Based
on a result from \cite{HI}\,, the constant $\sqrt{2}$ turns out to
be optimal in this case, and for the sequence
$\{e^{i2^nt}\}_{n=1}^\infty$\,. We also consider the case of a
sequence of Rademacher functions, and prove that the corresponding
noncommutative Khintchine inequality holds with constant $\sqrt{3}$
instead of $\sqrt{2}$\,, but we do not know yet whether this is
sharp.

In the second part of our paper we obtain an improvement of a recent
result of M. Junge (cf. \cite{Ju}) concerning a Khintchine-type
inequality for subspaces of $R{\oplus}_\infty C$ (the $l^\infty$-sum
of the row and column Hilbert spaces). Recall that
$R:=\text{Span}\{e_{1j}; j\geq 1\}$\,, respectively,
$C:=\text{Span}\{e_{j1}; j\geq 1\}$, where $e_{kl}$ is the element
in ${\mathcal B}(l_2)$ corresponding to the matrix with entries
equal to $1$ on the $(k, l)$ position, and $0$ elsewhere. This
Khintchine-type inequality is intimately connected with the question
of the existence of a completely isomorphic embedding of the
operator space $OH$, introduced by G. Pisier (see \cite{Pi6}), into
a noncommutative $L^1$-space, a problem that was resolved by Junge
in the remarkable paper \cite{Ju2}. In \cite{Ju} (see Section 8),
Junge improved this result, by showing that $OH$ cb-embeds into the
predual of a hyperfinite type III$_1$ factor.

In our new approach, we first observe that given a closed subspace
$H$ of $R\oplus C$\,, there is a self-adjoint operator $A\in
{\mathcal B}(H)$ satisfying $0\leq A\leq I$\,, where $I$ denotes the
identity operator on $H$\,, such that the operator space structure
on $H$ is given by
\begin{eqnarray}
\left\|\sum\limits_{i=1}^r x_i\otimes {\xi}_i
\right\|_{M_n({H})}&=&\max\left\{ \left\|\sum\limits_{i, j=1}^r
\langle (I-A){\xi}_i, {\xi}_j\rangle_H x_ix_j^*\right\|^{1/2}\,,
\left\|\sum\limits_{i, j=1}^r \langle A{\xi}_i, {\xi}_j\rangle_H
x_i^*x_j\right\|^{1/2} \right\}\,,
\end{eqnarray}
where $n, r$ are positive integers, $x_1\,, \ldots \,,x_r\in
M_n(\mathbb{C})$ and $\xi_1\,, \ldots \,,\xi_r\in H$\,.

As in Junge's approach from \cite{Ju}, we will use CAR algebra
methods. We consider the associated quasi-free state $\omega_A$ on
the CAR-algebra $\A=\A(H)$ built on the Hilbert space $H$\,, and
construct a linear map $F_A$ of $H^*$ into the predual $M_*$ of the
von Neumann algebra $M:=\overline{\pi_A(\A)}^{\text{sot}}$, which by
\cite{PS} is a hyperfinite factor. Here $\pi_A$ is the unital
$*$-homomorphism from the GNS representation associated to $(\A,
\omega_A)$. Note that $M_*$ can be considered as a subspace of
$\A^*$\,. Next we let $F_A$ be the transpose of the map
$E_A:\A\rightarrow H$ defined by
\begin{eqnarray}
\langle E_A(b), f\rangle_H&=&\omega_A(b{a(f)}^*+{a(f)}^*b)\,,
\quad\forall b\in \A\,, \forall f\in H\,,
\end{eqnarray}
where $f\mapsto a(f)$ is the map from $H$ to $\A=\A(H)$ in the
definition of the $CAR$-algebra (cf. \cite{D}). We then prove that
$F_A$ is a cb-isomorphism of $H^*$ onto its range, satisfying the
following estimates
\begin{equation}\label{eq102}
\frac1{\sqrt{2}}\|F_A(y)\|_{{M_n(\A)}^*}\leq \|y\|_{{M_n(H)}^*}\leq
\|F_A(y)\|_{{M_n(\A)}^*}\,, \quad\forall n\in \mathbb{N}\,, y\in
{M_n(H)}^*\,.
\end{equation}
We do so by first proving the dual version of the inequalities
(\ref{eq102}), namely we show that
\begin{equation}\label{eq1034}
\left\| \xi\right\|_{M_n(H)}\leq \left\|(\text{Id}_n\otimes
q_A)(\xi)\right\|_{M_n({{\A}/{\text{Ker}(E_A)}})}\leq
\sqrt{2}\left\|\xi\right\|_{M_n(H)}\,, \quad\forall n\in
\mathbb{N}\,, \xi\in M_n(H)\,.
\end{equation}
The estimate of the upper bound $\sqrt{2}$ in (\ref{eq1034})
(corresponding to the lower bound $\frac1{\sqrt{2}}$ in
(\ref{eq102})) is obtained by methods very similar to those we used
for the Pisier and Lust-Piquard noncommutative Khintchine
inequality. We then prove that both constants in (\ref{eq102}) are
sharp.

Note that if $P$ is the unique hyperfinite factor of type III$_1$
(cf. \cite{Ha3})\,, then the von Neumann algebra tensor product
$M\bar{\otimes} P$ is isomorphic to $P$\,, and therefore $F_A$ can
be considered as a completely bounded embedding of $H^*$ into the
predual $P_*$ of $P$\,, as well. It follows that every subspace of a
quotient of $(R\oplus C)^*$ is cb-isomorphic to a subspace of $P_*$
with cb-isomorphism constant $\leq \sqrt{2}$\,. In particular, due
to results of G. Pisier (cf. \cite{Pi2} Proposition A1)\,, the
operator Hilbert space $OH$ has this property (cf. Corollary
\ref{ohembed} in this paper). The question whether $OH$ embeds
completely isometrically into a noncommutative $L^1$-space remains
open.

In the case when the self-adjoint operator $A$ associated to the
subspace $H$ of $R\oplus C$ has pure point spectrum and
$\text{Ker}(A)=\text{Ker}(I-A)=0$\,, our construction of the map
$F_A: H^*\rightarrow M_*$  is very similar to Junge's construction
from \cite{Ju}\,. This can be seen by taking Lemma
\ref{lem677767776} into account.

We refer to the monographs \cite{ER, Pi1} for details on operator
spaces. We shall briefly recall some definitions that are relevant
for our paper. An operator space $V$ is a Banach space given
together with an isometric embedding $V\subset {\mathcal B}(H)$\,,
the algebra of bounded linear operators on a Hilbert space $H$\,.
For all $n\in \mathbb{N}$, this embedding determines a norm on
$M_n(V)$, the algebra of $n\times n$ matrices over $V$, induced by
the space $M_n({\mathcal B}(H))\cong {\mathcal B}(H^n)$\,. If $W$ is
a closed subspace of $V\,,$ then both $W$ and $V/W$ are operator
spaces; the matrix norms on $V/W$ are defined by
$M_n(V/W)=M_n(V)/{M_n(W)}$\,. The morphisms in the category of
operator spaces are {\em completely bounded maps}. Given a linear
map $\phi : V_0 \rightarrow V_1$ between two operator spaces $V_0$
and $V_1$\,, define ${\phi}_n : {M_n(V_0)} \rightarrow {M_n(V_1)}$
by ${\phi}_n([v_{ij}])= [\phi(v_{ij})]$\,, for all $
[v_{ij}]_{i,j=1}^n\in M_n(V_0)$\,. Let $\|\phi\|_{cb}:=\sup
\{\|{\phi}_n\|\,; n\in {\mathbb{N}} \,\}$\,. The map $\phi$ is
called {\em completely bounded} (for short, {\em cb}) if
$\|\phi\|_{cb} <\infty\,,$ and $\phi$ is called {\em completely
isometric} if all ${\phi}_n$ are isometries. A cb map $\phi$ which
is invertible with a cb inverse is called a {\em cb isomorphism}.
The space of all completely bounded maps from $V_0$ to $V_1$\,,
denoted by $\mathcal{C}\B(V_0, V_1)$\,, is an operator space with
matrix norms defined by $M_n(\mathcal{C}\B(V_0,
V_1))=\mathcal{C}\B(V_0, M_n(V_1))$\,. The dual of an operator space
$V$ is, again, an operator space $V^*=\mathcal{C}\B(V,
\mathbb{C})$\,.

\section{The Pisier and Lust-Piquard noncommutative Khintchine inequality}
\setcounter{equation}{0}
I. The complex Gaussian case

Let $\{\gamma_n\}_{n\geq 1}$ be a sequence of independent standard
complex-valued Gaussian random variables on some probability space
$(\Omega, \mathbb{P})$. Recall that a complex-valued random variable
on $(\Omega, \mathbb{P})$ is called Gaussian {\em standard} if it
has density $\frac1{\pi}e^{-|z|^2} d\text{Re} z \,d\text{Im} z\,.$
Equivalently, its real and imaginary parts are real-valued,
independent Gaussian random variables on $(\Omega, \mathbb{P})$\,,
each having mean 0 and variance $\frac12$\,. Therefore, for all
$n\geq 1$\,, $\mathbb{E}(\gamma_n)=0$ and
$\mathbb{E}(|\gamma_n|^2)=1$\,, where $\mathbb{E}$ denotes the usual
expectation of a random variable.

\begin{theorem}\label{th5}
Let $d$ and $n$ be positive integers, and consider $x_1, \ldots,
x_d\in M_n(\mathbb{C})$\,. Then the following inequalities hold
\begin{eqnarray}\label{eq80099}
{\frac1{\sqrt{2}}}|||\{x_i\}_{i=1}^d|||^*&\leq
&\left\|\sum\limits_{i=1}^d x_i\otimes \gamma_i\right\|_{L^1(\Omega;
S_1^n)}\leq |||\{x_i\}_{i=1}^d|||^*\,,
\end{eqnarray}
where $|||\{x_i\}_{i=1}^d|||^*$ is defined by (\ref{eq55565676})\,.
\end{theorem}

We will prove Theorem \ref{th5} by obtaining first its dual version,
namely,

\begin{prop}\label{dualcomplexgauss}
Let $d$ be a positive integer, and let $\{\gamma_i\}_{i=1}^d$ be a
sequence of independent standard complex-valued Gaussian random
variables on a probability space $(\Omega, \mathbb{P})$\,. For
$1\leq i\leq d$ define a map $\phi_i:L^\infty(\Omega)\rightarrow
\mathbb{C}$ by
\begin{eqnarray*}
\phi_i(f)&=&\int_\Omega f(\omega)\overline{\gamma_i}(\omega)
d\mathbb{P}(\omega)\,, \quad \forall \,f\in L^\infty(\Omega)\,,
\end{eqnarray*}
and let $E:L^\infty(\Omega)\rightarrow {\mathbb{C}}^d$ be defined by
\begin{eqnarray*}
E(f)&=&(\phi_1(f)\,, \ldots ,\phi_d(f))\,, \quad \forall \,f\in
L^\infty(\Omega)\,.
\end{eqnarray*}
Furthermore, let $q:L^\infty(\Omega)\rightarrow
{L^\infty(\Omega)}/{\text{Ker}(E)}$ denote the quotient map. Then,
for any positive integer $n$ and any $X\in M_n(L^\infty(\Omega))$\,,
\begin{eqnarray}\label{eq22334455}
\left|\left|\left|
\{x_i\}_{i=1}^d\right|\right|\right|_{M_n(\mathbb{C}^d)}&\leq
\|(\text{Id}_n\otimes
q)(X)\|_{M_n({L^\infty(\Omega)}/{\text{Ker}(E)})}\leq
\sqrt{2}\left|\left|\left|\{x_i\}_{i=1}^d\right|\right|\right|_{M_n(\mathbb{C}^d)}\,,
\end{eqnarray}
where $x_i=(\text{Id}_n\otimes \phi_i)(X)$\,, $\forall 1\leq i\leq
d$\,, and \begin{eqnarray} \label{eq99900099900} \left|\left|\left|
\{x_i\}_{i=1}^d\right|\right|\right|_{M_n(\mathbb{C}^d)}:&=&
\max\left\{\left\|\sum\limits_{i=1}^d x_i^*x_i\right\|^{\frac12}\,,
\left\|\sum\limits_{i=1}^d x_ix_i^*\right\|^{\frac12}\right\}\,.
\end{eqnarray}
\end{prop}

Note that ${\mathbb{C}}^d$ equipped with the sequence of matrix
norms $\{|||\cdot|||_{M_n(\mathbb{C})}\,, n\in \mathbb{N}\})$ is an
operator space.

\begin{proof}
Let $n\in \mathbb{N}$\,. We first prove the left hand side
inequality in (\ref{eq22334455}). For this, we need the following
\begin{lemma}\label{lem3131}
Let $X\in M_n(L^\infty(\Omega))$\,, and set
$x_i:=(\text{Id}_n\otimes \phi_i)(X)$\,, $\forall 1\leq i\leq d$\,.
Then
\begin{eqnarray}\label{eq22334456}
\|X\|_{M_n(L^\infty(\Omega))}&\geq
&\max\left\{\left\|\sum\limits_{i=1}^d x_i^*x_i\right\|^{\frac12}\,,
\left\|\sum\limits_{i=1}^d x_ix_i^*\right\|^{\frac12}\right\}\,.
\end{eqnarray}
\end{lemma}
\begin{proof}
Since $X\in M_n(\mathbb{C})\otimes L^\infty(\Omega)$
(algebraic tensor product), we can write
\[ X=\sum\limits_{k=1}^r y_k\otimes f_k\,, \]
for some $y_k\in M_n(\mathbb{C})$\,, $ f_k\in L^\infty(\Omega)$\,,
$1\leq k\leq r$\,.

Let $K=\text{Span}\{\gamma_1\,, \ldots \,,\gamma_d\,, f_1\,, \ldots
\,,f_r\}\subset L^2(\Omega)$\,. Choose an orthonormal basis
$\{g_i\}_{i=1}^s$ for $K$ such that
\begin{equation}\label{eq22334457}
g_i=\gamma_i\,, \quad 1\leq i\leq d\,.
\end{equation}
Then $X=\sum\limits_{i=1}^s z_i\otimes g_i$\,, for some $z_i\in
M_n(\mathbb{C})$\,, $1\leq i\leq s$\,. Note that for $1\leq i\leq
d$\,, we have \[ x_i=(\text{Id}_n\otimes
\phi_i)\left(\sum\limits_{j=1}^s z_j\otimes
g_j\right)=\sum\limits_{j=1}^sz_i\phi_i(g_j)=z_i\,, \] because by
(\ref{eq22334457}) it follows that $\phi_i(g_j)=\langle g_j,
\gamma_i\rangle_{L^2(\Omega)}=\langle g_j,
g_i\rangle_{L^2(\Omega)}=\delta_{ij}$\,, for all $1\leq j\leq s$\,.

Denote by $S(M_n(\mathbb{C}))$ the state space of
$M_n(\mathbb{C})$\,. Then, for $\omega\in S(M_n(\mathbb{C}))$\,,
\begin{eqnarray*}
\|X\|^2_{M_n(L^\infty(\Omega))}&\geq &(\omega\otimes
\mathbb{E})(X^*X)\\&=&(\omega\otimes
\mathbb{E})\left(\sum\limits_{i, j=1}^s z_i^*z_j\otimes
\bar{g_i}g_j\right)\\&=&\omega\left(\sum\limits_{i=1}^s
z_i^*z_i\right) \geq \omega\left(\sum\limits_{i=1}^d
z_i^*z_i\right)=\omega\left(\sum\limits_{i=1}^d x_i^*x_i\right)\,.
\end{eqnarray*}
Take supremum over all $\omega\in S(M_n(\mathbb{C}))$ to obtain
\begin{eqnarray}\label{eq22334458}
\|X\|^2_{M_n(L^\infty(\Omega))}&\geq &\left\|\sum\limits_{i=1}^d
x_i^*x_i\right\|\,.
\end{eqnarray}
Since
$\|X\|^2_{M_n(L^\infty(\Omega))}=\|XX^*\|_{M_n(L^\infty(\Omega))}$\,,
a similar argument shows that also
\begin{eqnarray}\label{eq22334459}
\|X\|^2_{M_n(L^\infty(\Omega))}&\geq &\left\|\sum\limits_{i=1}^d
x_ix_i^*\right\|\,.
\end{eqnarray}
This proves the lemma.
\end{proof}

\begin{rem}\label{rem89}\rm
 As a consequence of this lemma, we deduce
that for all $X\in M_n(L^\infty(\Omega))$ we have
\begin{eqnarray}\label{eq22334460}
\left|\left|\left|\{x_i\}_{i=1}^d\right|\right|\right|_{M_n(\mathbb{C}^d)}&\leq
& \|(\text{Id}_n\otimes
q)(X)\|_{M_n({L^\infty(\Omega)}/{\text{Ker}(E)})}\,,
\end{eqnarray}
i.e., the left hand side inequality in (\ref{eq22334455}) holds.
Indeed, for any $Y\in M_n(\text{Ker}(E))$ we infer by
(\ref{eq22334456}) that
\begin{equation*}
\|X+Y\|_{M_n(L^\infty(\Omega))}\,\geq \,
\left|\left|\left|(\text{Id}_n\otimes
E)(X+Y)\right|\right|\right|_{M_n(\mathbb{C}^d)}
=\left|\left|\left|(\text{Id}_n\otimes
E)(X)\right|\right|\right|_{M_n(\mathbb{C}^d)}
=\left|\left|\left|\{x_i\}_{i=1}^d\right|\right|\right|_{M_n(\mathbb{C}^d)}\,.
\end{equation*}
By taking infimum over all $Y\in M_n(\text{Ker}(E))$, inequality
(\ref{eq22334460}) follows by the definition of the quotient
operator space norm.
\end{rem}

It remains to prove the right hand side inequality in
(\ref{eq22334455})\,. For this, let $y_1\,, \ldots \,, y_d\in
M_n(\mathbb{C})$ and set
\[ Y:=\sum\limits_{i=1}^d y_i\otimes \gamma_i\in M_n(L^4(\Omega))\,. \] We will first compute
$(\text{Id}_n\otimes \mathbb{E})(Y^*Y)$\,, $(\text{Id}_n\otimes
\mathbb{E})(YY^*)$\,, $(\text{Id}_n\otimes \mathbb{E})((Y^*Y)^2)$
and $(\text{Id}_n\otimes \mathbb{E})((YY^*)^2)$\,.

Since $\mathbb{E}(\overline{\gamma_i}\gamma_j)=\delta_{ij}$\,, for
all $1\leq i, j\leq d$\,, we immediately get
\begin{eqnarray}\label{eq220022002}
(\text{Id}_n\otimes \mathbb{E})(Y^*Y)=\sum\limits_{i=1}^d
y_i^*y_i\,, && (\text{Id}_n\otimes
\mathbb{E})(YY^*)=\sum\limits_{i=1}^d y_iy_i^*\,.
\end{eqnarray}
It is easily checked that the vectors
$f_{ij}:=\overline{\gamma_i}\gamma_j-\delta_{ij}1$\,, $1\leq i,
j\leq d$\,, together with the constant function $1$ form an
orthonormal set with respect to the usual $L_2(\Omega)$-inner
product. We then obtain the expansion
\[ Y^*Y=\sum\limits_{i, j=1}^d y_i^*y_j\otimes
f_{ij}+\sum\limits_{i=1}^d y_i^*y_i\otimes 1\,, \] from which we
infer that
\begin{eqnarray}\label{eq200002}
(\text{Id}_n\otimes \mathbb{E})((Y^*Y)^2)&=&\sum\limits_{i=1}^d
y_i^*\left(\sum\limits_{j=1}^d y_jy_j^*\right)y_i
+\left(\sum\limits_{i=1}^d y_i^* y_i\right)^2\,.
\end{eqnarray}
A similar argument shows that
\begin{eqnarray}\label{eq200003}
(\text{Id}_n\otimes \mathbb{E})((YY^*)^2)&=& \sum\limits_{i=1}^d
y_i\left(\sum\limits_{j=1}^d
y_j^*y_j\right)y_i^*+\left(\sum\limits_{i=1}^d y_iy_i^*\right)^2\,.
\end{eqnarray}
By (\ref{eq200002})\,, (\ref{eq200003}) and (\ref{eq220022002}) we
then obtain the following inequalities
\begin{eqnarray}
(\text{Id}_n\otimes \mathbb{E})((Y^*Y)^2)&\leq &
\left(\left\|\sum\limits_{i=1}^d
y_i^*y_i\right\|+\left\|\sum\limits_{i=1}^d
y_iy_i^*\right\|\right)(\text{Id}_n\otimes
\mathbb{E})(Y^*Y)\,,\label{eq34343446}\\
(\text{Id}_n\otimes \mathbb{E})((YY^*)^2)&\leq &
\left(\left\|\sum\limits_{i=1}^d
y_i^*y_i\right\|+\left\|\sum\limits_{i=1}^d
y_iy_i^*\right\|\right)(\text{Id}_n\otimes
\mathbb{E})(YY^*)\,.\label{eq34343445}
\end{eqnarray}
The crucial point in proving the right hand side inequality in
(\ref{eq22334455}) is to show the following
\begin{lemma}\label{lem3232}
Let $x_1\,, \ldots \,,x_d\in M_n(\mathbb{C})$\,. Then there exists
$X\in M_n(L^\infty(\Omega))$ such that
\[ (\text{Id}_n\otimes E)(X)=\sum \limits_{i=1}^d x_i\otimes e_i\,, \]
where $\{e_i\}_{1\leq i\leq d}$ is the canonical unit vector basis in $\mathbb{C}^d$\,, and
\begin{eqnarray*}
\|X\|_{M_n(L^\infty(\Omega))}&\leq & \sqrt{2}
\max\left\{\left\|\sum\limits_{i=1}^d x_i^*x_i\right\|^{\frac12}\,,
\left\|\sum\limits_{i=1}^d x_ix_i^*\right\|^{\frac12}\right\}\,.
\end{eqnarray*}
\end{lemma}

We first prove the following lemma:
\begin{lemma}\label{lem4567}
If $y_1\,, \ldots \,,y_d\in M_n(\mathbb{C})$ and
\begin{eqnarray}\label{eq6677889988}
\max\left\{\left\|\sum\limits_{i=1}^d y_i^*y_i\right\|^{\frac12}\,,
\left\|\sum\limits_{i=1}^d y_iy_i^*\right\|^{\frac12}\right\}&=&1\,,
\end{eqnarray} then there exists $Z\in M_n(L^\infty(\Omega))$ such
that
\begin{eqnarray*}
\|Z\|_{M_n(L^\infty(\Omega))}&\leq & \frac1{\sqrt{2}}
\end{eqnarray*}
and, moreover, when $z_1\,, \ldots \,,z_d$ are defined by
$(\text{Id}_n\otimes E)(Z)=\sum\limits_{i=1}^d z_i\otimes e_i$\,,
then
\begin{eqnarray*}
\max\left\{\left\|\sum\limits_{i=1}^d
(y_i-z_i)^*(y_i-z_i)\right\|^{\frac12}\,, \left\|\sum\limits_{i=1}^d
(y_i-z_i)(y_i-z_i)^*\right\|^{\frac12}\right\}&\leq & \frac12\,.
\end{eqnarray*}
\end{lemma}

\begin{proof}
Set
\begin{eqnarray*}
Y&=& \sum\limits_{i=1}^d y_i\otimes \gamma_i\in M_n(L^4(\Omega))\,.
\end{eqnarray*}
Let $\widetilde{E}: L^4(\Omega)\rightarrow \mathbb{C}^d$ denote the
natural extension of $E$ to $L^4(\Omega)$\,. Then
\[ (\text{Id}_n\otimes \widetilde{E})(Y)=\sum\limits_{i=1}^d y_i\otimes \widetilde{E}(\gamma_i)=\sum\limits_{i=1}^d y_i\otimes e_i\,. \]
Now let $C>0$ and define $F_C:\mathbb{R}\rightarrow \mathbb{R}$ by
\begin{eqnarray}\label{eq7070709822}
F_C(t)&=&\left\{\begin{array}{lll}
                                 -C& \,\mbox{if} \;\; \;t< -C\\
                                       t & \,\mbox{if} \;\;\;-C\leq t\leq C\\
                                       C & \,\mbox{if}\;\;\;t>C\\
                                \end{array}
                        \right.
\end{eqnarray}
Use functional calculus to define $Z\in M_n(L^\infty(\Omega))$ by
\begin{eqnarray}\label{eq34345656560}
\left(\begin{array}
[c]{cc}%
0 & Z^*\\
Z & 0
\end{array}
\right)&=&F_C\left(\begin{array}
[c]{cc}%
0 & Y^*\\
Y & 0
\end{array}
\right)\,.
\end{eqnarray}
Note that this implies that $\|Z\|_{M_n(L^\infty(\Omega))}\leq C$\,.
Further, set
\begin{eqnarray*}
G_C(t)&=& t-F_C(t)\,, \quad \forall t\in \mathbb{R}\,.
\end{eqnarray*}
We then have
\begin{eqnarray*}
\left(\begin{array}
[c]{cc}%
0 & (Y-Z)^*\\
(Y-Z) & 0
\end{array}
\right)=G_C\left(\begin{array}
[c]{cc}%
0 & Y^*\\
Y & 0
\end{array}
\right)
\end{eqnarray*}
and thus
\begin{eqnarray}\label{eq567890}
\left(\begin{array}
[c]{cc}%
(Y-Z)^*(Y-Z) & 0\\
0 & (Y-Z)(Y-Z)^*
\end{array}
\right)=\left(G_C\left(\begin{array}
[c]{cc}%
0 & Y^*\\
Y & 0
\end{array}
\right)\right)^2\,.
\end{eqnarray}
A simple calculation shows that
\begin{eqnarray}\label{eq567891}
|G_C(t)|&\leq & \frac1{4C}t^2\,, \quad \forall t\in \mathbb{R}\,.
\end{eqnarray}
By functional calculus it follows that
\begin{eqnarray*}\label{eq1000000}
\left(G_C\left(\begin{array}
[c]{cc}%
0 & Y^*\\
Y & 0
\end{array}
\right)\right)^2&\leq &\frac1{16C^2}{\left(\begin{array}
[c]{cc}%
0 & Y^*\\
Y & 0
\end{array}
\right)}^4 = \frac1{16C^2}{\left(\begin{array}
[c]{cc}%
(Y^*Y)^2 & 0\\
0 & (YY^*)^2
\end{array}
\right)}\,.\nonumber
\end{eqnarray*}
Hence, by (\ref{eq567890}) and (\ref{eq567891}) we infer that
\begin{eqnarray}
(Y-Z)^*(Y-Z)&\leq & \frac1{16C^2}(Y^*Y)^2\,,\label{eq9988998871}\\
(Y-Z)(Y-Z)^*&\leq & \frac1{16C^2}(YY^*)^2\,.\label{eq9988998881}
\end{eqnarray}
By letting $z_i=(\text{Id}_n\otimes \phi_i)(Z)$\,, $1\leq i\leq d$\,, we then have
\[ (\text{Id}_n\otimes E)(Z)=\sum\limits_{i=1}^d z_i\otimes e_i\,,
\]
and hence $(\text{Id}_n\otimes
\widetilde{E})(Y-Z)=\sum\limits_{i=1}^d (y_i-z_i)\otimes e_i$\,.

By (\ref{eq9988998871})\,, (\ref{eq34343446})  and (\ref{eq6677889988}) we then obtain the estimates
\begin{eqnarray*}
\sum\limits_{i=1}^d (y_i-z_i)^*(y_i-z_i)&\leq & (\text{Id}_n\otimes \mathbb{E})((Y-Z)^*(Y-Z))\\
&\leq & \frac1{16C^2}(\text{Id}_n\otimes \mathbb{E})((Y^*Y)^2)\\
&\leq & \frac1{16C^2}\left(\left\|\sum\limits_{i=1}^d y_i^*y_i\right\|+\left\|\sum\limits_{i=1}^d y_iy_i^*\right\|\right)                  (\text{Id}_n\otimes \mathbb{E})(Y^*Y)\\
&\leq & \frac{2}{16C^2}(\text{Id}_n\otimes \mathbb{E})(Y^*Y)\\&=&
\frac1{8C^2}\sum\limits_{i=1}^d y_i^*y_i\,.
\end{eqnarray*}
It follows that
\[
\left\|\sum\limits_{i=1}^d (y_i-z_i)^*(y_i-z_i)\right\|\leq \frac1{8C^2}\left\|\sum\limits_{i=1}^d y_i^*y_i\right\|
\leq  \frac1{8C^2}\,. \]
Similarly, we also get
\begin{eqnarray*}
\left\|\sum\limits_{i=1}^d (y_i-z_i)(y_i-z_i)^*\right\|
&\leq & \frac1{8C^2}\,.
\end{eqnarray*}
Hence,
\begin{eqnarray*}
\max\left\{\left\|\sum\limits_{i=1}^d
(y_i-z_i)^*(y_i-z_i)\right\|^{\frac12}\,, \left\|\sum\limits_{i=1}^d
(y_i-z_i)(y_i-z_i)^*\right\|^{\frac12}\right\} &\leq &
\frac1{\sqrt{8}C}\,.
\end{eqnarray*}
Now take $C=\frac1{\sqrt{2}}$ to get the conclusion.
\end{proof}
We also need the following result:
\begin{lemma}\label{lem25}
Let $V$ and $W$ be Banach spaces. Consider $T:V\rightarrow W$ a
bounded linear map. Further, let $\phi:W\rightarrow V$ be a
non-linear map such that, for some $C> 0$ and some $0< \delta< 1$\,,
we have
\begin{eqnarray}
\|\phi(w)\|&\leq & C\|w\|\,,\label{eq67678873}\\
\|w-(T\circ \phi)(w)\|&\leq &\delta\|w\|\,, \quad \forall w\in W\,.\label{eq67678874}
\end{eqnarray}
Then there exists a non-linear map $\psi:W\rightarrow V$ such that
$T\circ \psi=\text{Id}_W$ and, moreover,
\begin{eqnarray*}
\|\psi(w)\|&\leq &\frac{C}{1-\delta}\|w\|\,, \quad\forall w\in W\,.
\end{eqnarray*}
\end{lemma}

\begin{proof} Let $w\in W$\,. Set $w_0=w$ and define recursively
\begin{eqnarray*}
w_n&=&w_{n-1}-(T\circ \phi)(w_{n-1})\,, \quad \forall n\geq 1\,.
\end{eqnarray*}
Then, by (\ref{eq67678874}) we have for all $n\geq 0$
\begin{equation}\label{eq69}
\|w_{n+1}\|\leq \delta\|w_n\|\leq \ldots \leq \delta^{n+1}\|w_0\|\,.
\end{equation}
Also, we deduce that
\[
w=w_0=(T\circ \phi)(w_0)+w_1 =(T\circ \phi)(w_0)+(T\circ
\phi)(w_1)+w_2 =\ldots =\sum\limits_{j=0}^n (T\circ
\phi)(w_j)+w_{n+1}\,, \quad n\geq 0\,. \] By (\ref{eq69}) it follows
that $w_n\rightarrow 0$ as $n\rightarrow \infty$\,, and therefore
\begin{equation}\label{eq70}
w=\sum\limits_{j=0}^\infty (T\circ
\phi)(w_j)=T\left(\sum\limits_{j=0}^\infty \phi(w_j)\right)\,.
\end{equation}
Define
\[ \psi(w):= \sum\limits_{j=0}^\infty \phi(w_j)\,, \quad \forall
w\in W\,. \] By (\ref{eq70}) it follows that $T(\psi(w))=w$\,, for
all $w\in W$\,. Moreover, by (\ref{eq67678873}) and (\ref{eq69}) we
obtain that
\[ \|\psi(w)\|=\left\|\sum\limits_{j=0}^\infty
\phi(w_j)\right\|\leq C\sum\limits_{j=0}^\infty \|w_j\|\leq
\frac{C}{1-\delta}\|w\|\,, \quad \forall w\in W\,, \] which completes the
proof.
\end{proof}

Now we are ready to prove Lemma \ref{lem3232}\,. Indeed, Lemma \ref{lem4567} shows
that if
\[ y=\sum\limits_{i=1}^d y_i\otimes e_i\in M_n(\mathbb{C}^d) \]
satisfies $|||y|||_{M_n(\mathbb{C}^d)}=1$\,, then there exists $Z\in
M_n(L^\infty(\Omega))$ so that $\|Z\|_{M_n(L^\infty(\Omega))}\leq
\frac1{\sqrt{2}}$ and $|||(\text{Id}_n \otimes
E)(Z)-y|||_{M_n(\mathbb{C}^d)}\leq \frac12$\,. By homogeneity we
infer that for all $y\in M_n(\mathbb{C}^d)$ there exists $Z\in
M_n(L^\infty(\Omega))$ so that $\|Z\|_{M_n(L^\infty(\Omega))}\leq
\frac1{\sqrt{2}}|||y|||_{M_n(\mathbb{C}^d)}$\,, and, moreover,
$|||(\text{Id}_n \otimes E)(Z)-y|||_{M_n(\mathbb{C}^d)}\leq
{\frac12}|||y|||_{M_n(\mathbb{C}^d)}$\,. Apply now Lemma \ref{lem25}
with $V=M_n(L^\infty(\Omega))$\,, $W=M_n(\mathbb{C}^d)$\,,
$T=\text{Id}_n\otimes E$\,, the map $\phi:W\rightarrow V$ be defined
by $\phi(y)=Z$\,, $\forall y\in W$\,, $C=\frac1{\sqrt{2}}$ and
$\delta=\frac12$\,. We deduce that for all $x_1\,, \ldots \,, x_d\in
M_n(\mathbb{C})$\,, there exists $X\in M_n(L^\infty(\Omega))$ such
that $(\text{Id}_n\otimes E)(X)=\sum\limits_{i=1}^d x_i\otimes e_i$
and
\begin{equation}\label{eq898989675}\|X\|_{M_n(L^\infty(\Omega))}\leq
{\frac{C}{1-\delta}}\left|\left|\left|\sum\limits_{i=1}^d x_i\otimes
e_i\right|\right|\right|_{M_n(\mathbb{C}^d)}=\sqrt{2}|||\{x_i\}_{i=1}^d|||_{M_n(\mathbb{C}^d)}\,,
\end{equation} which completes the proof of Lemma \ref{lem3232}.
Note that, since the norm on
${M_n\left({L^\infty(\Omega)}/{\text{Ker}(E)}\right)}$ is the
quotient space norm on the space
${M_n(L^\infty(\Omega))}/{M_n(\text{Ker}(E))}$\,, it follows by
(\ref{eq898989675}) that
\[ \left\|(\text{Id}_n\otimes q)(X)\right\|_{M_n\left({L^\infty(\Omega)}/{\text{Ker}(E)}\right)}\leq \sqrt{2}
|||\{x_i\}_{i=1}^d|||_{M_n(\mathbb{C}^d)}\,. \] Therefore, by Lemmas
\ref{lem3131} and \ref{lem3232} and Remark \ref{rem89}\,, there is a
linear bijection
$\widehat{E}:{L^\infty(\Omega)}/{\text{Ker}(E)}\rightarrow
\mathbb{C}^d$ such that
\begin{eqnarray*}
\widehat{E}(q(s_i))&=& e_i\,, \quad \forall 1\leq i\leq d\,,
\end{eqnarray*}
where $s_i=\sqrt{\frac{4}{\pi}} \text{sgn}(\gamma_i)$\,, for $1\leq
i\leq d$\,. Note that $s_i\in L^\infty(\Omega)$ and
$\mathbb{E}(s_i\overline{\gamma_i})=\delta_{ij}$\,, $\forall 1\leq
j\leq d$\,, so that $E(s_i)=e_i$\,, $\forall 1\leq i\leq d$\,. For
every positive integer $n$\,, the following diagram is commutative,
$$
\xymatrix{
 {M_n(L^\infty(\Omega))}\ar@{->}^{\text{Id}_n\otimes E}[rr]
 \ar@{->}_{\text{Id}_n\otimes q}[dr]
 & & {M_n(\mathbb{C}^d)}\\
& {M_n({L^\infty(\Omega)}/{\text{Ker}(E)})}
\ar@{->}[ur]_{\text{Id}_n\otimes \widehat{E}}}
$$
and moreover, the inequalities (\ref{eq22334455}) hold. The proof of
Proposition \ref {dualcomplexgauss} is now complete.
\end{proof}

\begin{rem}\label{rem789}\rm
We should mention that, by the same proof with only minor
modifications, Theorem \ref{th5} remains valid if we replace the
sequence $\{\gamma_n\}_{n\geq 1}$ of independent standard complex
Gaussian random variables by a sequence $\{s_n\}_{n\geq 1}$ of
independent Steinhauss random variables (that is, a sequence of
independent random variables which are uniformly distributed over
the unit circle), or by the sequence $\{e_n\}_{n\geq 1}$ given by
$e_n(t)= e^{i{2^n}t}$\,, $0\leq t\leq 2\pi$\,. Indeed, the only
essential change in the proof is that the formulas (\ref{eq200002})
and (\ref{eq200003}) must be modified, because in the case of the
sequences $\{s_n\}_{n\geq 1}$ and $\{e_n\}_{n\geq 1}$ we still have
that $\{\bar{s}_is_j; 1\leq i, j\leq d\}\cup \{1\}$ and,
respectively, $\{\bar{e}_ie_j; 1\leq i, j\leq d\}\cup \{1\}$ form
orthonormal sets, but in contrast to the case of the Gaussians
$\{\gamma_n\}_{n\geq 1}$\,, one has
\[ \bar{s}_js_j=\bar{e}_je_j=1\,, \quad j\geq 1\,.
\]
Therefore, the diagonal terms (corresponding to $i=j$) in the right
hand sides of (\ref{eq200002}) and (\ref{eq200003}) should be
removed from the double sums. However, since the diagonal terms are
all positive, it follows that (\ref{eq34343446}) and
(\ref{eq34343445}) remain valid in the case of the sequences
$\{s_n\}_{n\geq 1}$ and $\{e_n\}_{n\geq 1}$\,, as well.
\end{rem}

We now discuss estimates for best constants in the noncommutative
Khintchine inequalities ($p=1$).

\begin{theorem}\label{th8998}
Denote by $c_1$\,, $c_2$ the best constants in the inequalities
\begin{eqnarray}\label{eq3322332244}
{c_1}|||\{x_i\}_{i=1}^d|||^*&\leq &\left\|\sum\limits_{i=1}^d
x_i\otimes \gamma_i\right\|_{L^1(\Omega; S_1^n)}\leq
{c_2}|||\{x_i\}_{i=1}^d|||^*\,,
\end{eqnarray}
where $d$ and $n$ are positive integers, $x_1\,, \ldots \,,x_d\in
M_n(\mathbb{C})$, and $\{\gamma_i\}_{i=1}^d$ is a sequence of
independent standard complex-valued Gaussian random variables on a
probability space $(\Omega, \mathbb{P})$\,. Then
\begin{equation*}
c_1=\frac1{\sqrt{2}}\,, \qquad c_2=1\,.
\end{equation*}
\end{theorem}

\begin{proof}
Let $m$ be a positive integer. Let $d=2m+1$ and set
$n={{2m+1}\choose {m}}$\,. Then, by Theorem 1.1 in \cite{HI}, there
exist partial isometries $a_1\,, \ldots \,,a_d\in {\mathcal
B}(H)$\,, where $H$ is a Hilbert space of $\text{dim}(H)=n$\,, such
that
\begin{eqnarray}\label{eq3434322}
\tau(a_i^*a_i)&=&\frac{m+1}{2m+1}\,, \quad \forall 1\leq i\leq d\,,
\end{eqnarray}
where $\tau$ denotes the normalized trace on ${\mathcal B}(H)$\,, satisfying, moreover,
\begin{equation}\label{eq5656776783}
\sum\limits_{i=1}^d a_i^*a_i=\sum\limits_{i=1}^d a_ia_i^*=(m+1)I\,,
\end{equation}
where $I$ denotes the identity operator on $H$\,. First, we claim
that
\begin{eqnarray}
|||\{a_i\}_{i=1}^d|||_{M_n(\mathbb{C}^d)}&=&\sqrt{m+1}\,,\label{eq34342234}\\
|||\{a_i\}_{i=1}^d|||^*&=&n\sqrt{m+1}\,.\label{eq34342235}
\end{eqnarray}
Indeed, (\ref{eq34342234}) follows immediately from the definition
of the norm $|||\cdot|||_{M_n(\mathbb{C}^d)}$ and relation
(\ref{eq5656776783}), while the equation (\ref{eq34342235}) follows
from the following estimates
\begin{eqnarray*}
|||\{a_i\}_{i=1}^d|||^*&=&\sup\left\{\left|\text{Tr}\left(\sum\limits_{i=1}^d a_ib_i\right)\right|; |||\{b_i\}_{i=1}^d|||_{M_n(\mathbb{C}^d)}\leq 1\right\}\\
&\geq & \left|\text{Tr}\left(\sum\limits_{i=1}^d a_i\left(\frac{a_i^*}{\sqrt{m+1}}\right)\right)\right|\\
&=& \frac1{\sqrt{m+1}}\text{Tr}\left(\sum\limits_{i=1}^d a_ia_i^*\right)\\
&=& \frac1{\sqrt{m+1}}\text{Tr}((m+1)I)=n\sqrt{m+1}\,,
\end{eqnarray*}
respectively,
\begin{equation*}
|||\{a_i\}_{i=1}^d|||^*\leq \text{Tr}\left(\left(\sum\limits_{i=1}^d
a_i^*a_i\right)^{1/2}\right) =\text{Tr}(\sqrt{m+1}I)=n\sqrt{m+1}\,.
\end{equation*}
It was proved in \cite{HI} that $a_1\,, \ldots \,,a_d$ have the
additional property that $\forall \beta_1\,, \ldots \,,\beta_d\in
\mathbb{C}$ with $\sum\limits_{i=1}^d |\beta_i|^2=1$\,, the operator
$y:=\sum\limits_{i=1}^d \beta_ia_i\in {\mathcal B}(H)$ is also a
partial isometry with $\tau(y^*y)=\frac{m+1}{2m+1}$\,. This implies
that for all $\omega\in \Omega$\,, the operator
\[ y_\omega:=\sum\limits_{i=1}^d \frac{\gamma_i(\omega)}{\left(\sum\limits_{i=1}^d |\gamma_i(\omega)|^2\right)^{\frac12}}\,a_i\in {\mathcal B}(H) \]
is a partial isometry with $\tau(y_\omega^*y_\omega)=\frac{m+1}{2m+1}$\,, and we deduce that
\begin{eqnarray}\label{eq5566776677881}
c_1|||\{a_i\}_{i=1}^d |||^*&\leq &\int_\Omega \left\|\sum\limits_{i=1}^d \gamma_i(\omega)a_i\right\|_{L^1(M_n(\mathbb{C})\,, \text{Tr})} d\mathbb{P}(\omega)\\
&=& \int_\Omega \left(\sum\limits_{i=1}^d |\gamma_i(\omega)|^2\right)^\frac12\|y_\omega\|_{L^1(M_n(\mathbb{C})\,, \text{Tr})}d\mathbb{P}(\omega)\nonumber\\
&=&n\cdot \frac{m+1}{2m+1}\int_\Omega \left(\sum\limits_{i=1}^d
|\gamma_i(\omega)|^2\right)^\frac12 d\mathbb{P}(\omega)\,,\nonumber
\end{eqnarray}
wherein we have used the fact that $y_\omega^*y_\omega$ is a
projection satisfying $\tau(|y_\omega|)=\frac{m+1}{2m+1}$\,, for all
$\omega\in \Omega$\,. A standard computation yields the formula
\begin{eqnarray}\label{eq5564644}
\int_\Omega \left(\sum\limits_{i=1}^d |\gamma_i(\omega)|^2\right)^\frac12 d\mathbb{P}(\omega)&=& \frac{\Gamma\left(d+\frac{1}{2}\right)}{\Gamma(d)}\,.
\end{eqnarray}
Indeed, since the distribution of $|\gamma_i|^2$ is $\Gamma(1,
1)$\,, $1\leq i\leq d$\,, it follows by independence that the
distribution of $\sum\limits_{i=1}^d |\gamma_i|^2$ is $\Gamma(d,
1)$, whose density is $\frac1{\Gamma(d)} x^{d-1}e^{-x}$\,,  $x>
0$\,. Since $\int_0 ^\infty x^{\frac12}x^{d-1}e^{-x}
dx=\Gamma(d+\frac12)$\,, formula (\ref{eq5564644}) follows.
Combining now (\ref{eq5566776677881}) with (\ref{eq34342235}) and
(\ref{eq5564644}) we deduce that
\begin{equation}\label{eq23322332111}
c_1\leq  \frac1{\sqrt{m+1}}\left(\frac{m+1}{2m+1}\right)\frac{\Gamma\left(d+\frac{1}{2}\right)}{\Gamma(d)}
\leq \frac{\sqrt{m+1}}{2m+1} \sqrt{2m+1}=\frac{\sqrt{m+1}}{\sqrt{2m+1}}\,,
\end{equation}
wherein we have used the inequality
\[ \Gamma\left(k+\frac12\right)< \sqrt{k} \Gamma(k)\,, \quad \forall k\in \mathbb{N}\,, \]
applied for $k=d=2m+1$\,. Since $m$ was arbitrarily chosen and
$\lim\limits_{m\rightarrow
\infty}\frac{\sqrt{m+1}}{\sqrt{2m+1}}=\frac1{\sqrt{2}}$\,, we deduce
by (\ref{eq23322332111}) that $c_1\leq \frac1{\sqrt{2}}$\,. By
Theorem  \ref{th5} we know that $c_1\geq \frac1{\sqrt{2}}$\,, hence
we conclude that $c_1=\frac1{\sqrt{2}}$\,.

To estimate $c_2$\,, let $d$ be a positive integer. Set $n=d$\,. For
all $1\leq i\leq d$, set $x_i:=e_{i1}\in M_d(\mathbb{C})$\,. We then
have
\[ \left\|\sum\limits_{i=1}^d
x_i\otimes \gamma_i\right\|_{L^1(\Omega; S_1^n)}=\int_\Omega
\left\|\sum\limits_{i=1}^d
\gamma_i(\omega)x_i\right\|_{L^1(M_n(\mathbb{C})\,, \text{Tr})}
d\mathbb{P}(\omega)=\int_\Omega \left(\sum\limits_{i=1}^d
|\gamma_i(\omega)|^2\right)^\frac12 d\mathbb{P}(\omega)\,. \] Note
also that
\[ |||\{x_i\}_{i=1}^d|||^*\geq
\text{Tr}\left(\left(\sum\limits_{i=1}^d
x_i^*x_i\right)^{\frac12}\right)=\text{Tr}(\sqrt{d}\,\,e_{11})=1\,.
\] Then, using (\ref{eq5564644}), together with the fact that
$\lim\limits_{d\rightarrow \infty} \frac{1}{\sqrt{d}}
\frac{\Gamma\left(d+\frac12\right)}{\Gamma(d)}=1$\,, we infer by
(\ref{eq3322332244}) that $c_2\geq 1$\,. Since by Theorem \ref{th5}
we get $c_2\leq 1$\,, we conclude that $c_2=1$\,, and the proof is
complete.
\end{proof}

\begin{rem}\label{rem78978}\rm
If we replace the sequence of independent standard complex-valued
Gaussian random variables $\{\gamma_n\}_{n\geq 1}$ by a sequence of
independent Steinhauss random variables $\{s_n\}_{n\geq 1}$ or by
the sequence $\{e^{i2^nt}\}_{n\geq 1}$\,, and denote by $c_1\,, c_2$
the best constants in the corresponding inequalities
(\ref{eq3322332244}), the same argument will give
$c_1=\frac1{\sqrt{2}}$\,. Also, $c_2=1$ in both cases, as a
consequence of Remark \ref{rem789} and the fact that
$\|s_1\|_{L^1(\mathbb{T})}=1=\|e^{i2t}\|_{L^1(\mathbb{T})}$\,, where
$\mathbb{T}$ is the unit circle with normalized Lebesgue measure
${dt}/{2\pi}$\,.
\end{rem}

\noindent
II. The Rademacher case\\

Let $\{r_n\}_{n\geq 1}$ be a sequence of Rademacher functions on
$[0, 1]$\,. Probabilistically, one can think of $\{r_n\}_{n\geq 1}$
as being a sequence of independent, identically distributed random
variables on $[0, 1]$, each taking value 1 with probability
$\frac12$, respectively, value $-1$ with probability $\frac12$\,. It
is easily seen that $\mathbb{E}(r_n)=0$ and
$\mathbb{E}(r_nr_m)=\delta_{nm}$\,, for all $n, m\in \mathbb{N}$\,.

\begin{theorem}\label{th4}
Let $d$ and $n$ be positive integers and consider $x_1, \ldots, x_d\in M_n(\mathbb{C})$\,.
Then the following inequalities hold
\begin{eqnarray}\label{eq8009}
{\frac1{\sqrt{3}}}|||\{x_i\}_{i=1}^d|||^*&\leq
&\left\|\sum\limits_{i=1}^d x_i\otimes r_i\right\|_{L^1([0, 1];
S_1^n)}\leq |||\{x_i\}_{i=1}^d|||^*\,.
\end{eqnarray}
\end{theorem}

As in the case of complex Gaussian random variables, we prove the
dual version of Theorem \ref{th4}, namely,

\begin{prop}\label{dualradem}
Let $d$ be a positive integer, and let $\{r_i\}_{1\leq i\leq
d}$ be a sequence of Rademacher functions
on $[0,1]$\,. For $1\leq i\leq d$
define $\phi_i:L^\infty([0,1])\rightarrow \mathbb{C}$ by
\begin{eqnarray*}
\phi_i(f)&=&\int_0^1 f(t)r_i(t)
dt\,, \quad \forall \,f\in L^\infty([0,1])\,,
\end{eqnarray*}
and let $E:L^\infty([0,1])\rightarrow {\mathbb{C}}^d$ be defined by
\begin{eqnarray*}
E(f)&=&(\phi_1(f)\,, \ldots ,\phi_d(f))\,, \quad \forall \,f\in
L^\infty([0,1])\,.
\end{eqnarray*}
Furthermore, let $q:L^\infty([0,1])\rightarrow
{L^\infty([0,1])}/{\text{Ker}(E)}$ denote the quotient map. Then,
for any positive integer $n$ and any $X\in M_n(L^\infty([0,1]))$\,,
\begin{eqnarray}\label{eq2233445577}
\left|\left|\left|
\{x_i\}_{i=1}^d\right|\right|\right|_{M_n(\mathbb{C}^d)}&\leq
\|(\text{Id}_n\otimes
q)(X)\|_{M_n({L^\infty([0,1])}/{\text{Ker}(E)})}\leq
\sqrt{3}\left|\left|\left|\{x_i\}_{i=1}^d\right|\right|\right|_{M_n(\mathbb{C}^d)}\,,
\end{eqnarray}
where $x_i=(\text{Id}_n\otimes \phi_i)(X)$\,, $\forall 1\leq i\leq
d$\,.
\end{prop}

\begin{proof}
Let $n$ be a positive integer. The proof of the left hand side
inequality in (\ref{eq2233445577}) is the same as in the complex
Gaussian case. For the right hand side inequality we follow the same
argument, but with appropriate modifications, which we indicate
below.

Let $y_1\,, \ldots\,, y_d\in M_n(\mathbb{C})$ and set
\[ Y:=\sum\limits_{i=1}^d y_i\otimes r_i\in M_n(L^\infty([0, 1])\,. \]
As before we will estimate $(\text{Id}_n\otimes
\mathbb{E})(Y^*Y)$\,, $(\text{Id}_n\otimes \mathbb{E})(YY^*)$\,,
$(\text{Id}_n\otimes \mathbb{E})((Y^*Y)^2)$ and $(\text{Id}_n\otimes
\mathbb{E})((YY^*)^2)$\,.

First note that $Y^*Y=\sum\limits_{i, j=1}^{d} y_i^*y_j\otimes r_ir_j$ and, respectively, $YY^*=\sum\limits_{i, j=1}^{d} y_iy_j^*\otimes r_ir_j$ to conclude that
\begin{eqnarray}
(\text{Id}_n\otimes \mathbb{E})(Y^*Y)&=& \sum\limits_{i=1}^{d} y_i^*y_i\label{eq11119}\\
(\text{Id}_n\otimes \mathbb{E})(YY^*)&=& \sum\limits_{i=1}^{d}
y_iy_i^*\,.\label{eq11120}
\end{eqnarray}
Furthermore, note that $(\text{Id}_n\otimes \mathbb{E})((Y^*Y)^2)=
\sum\limits_{i, j, k, l=1}^d y_i^*y_jy_k^*y_l
\mathbb{E}(r_ir_jr_kr_l)$\,. Since $\mathbb{E}(r_ir_jr_kr_l)\in \{0,
1\}$ with $\mathbb{E}(r_ir_jr_kr_l)=1$ if and only if $i=j=k=l$, or
$i=j\ne k=l$, or $i=k\ne j=l$, or $i=l\ne j=k$, it then follows that
\begin{eqnarray}\label{eq3333444}
(\text{Id}_n\otimes \mathbb{E})((Y^*Y)^2)&=& \sum\limits_{i=1}^d y_i^*y_iy_i^*y_i + \sum\limits_{i\ne j}y_i^*y_iy_j^*y_j + \sum\limits_{i\ne j} y_i^*y_jy_i^*y_j + \sum\limits_{i\ne j} y_i^*y_jy_j^*y_i\\
&=& \sum\limits_{i, j=1}^d y_i^*y_iy_j^*y_j + \sum\limits_{i\ne j} y_i^*y_jy_i^*y_j + \sum\limits_{i\ne j} y_i^*y_jy_j^*y_i\,.\nonumber
\end{eqnarray}
Note that $\sum\limits_{i, j=1}^d
y_i^*y_iy_j^*y_j=\left(\sum\limits_{i=1}^d y_i^*y_i\right)^2$\,.
Further, we have
\begin{equation*}
\sum\limits_{i\ne j}y_i^*y_jy_i^*y_j=\sum\limits_{i<
j}y_i^*y_jy_i^*y_j+ \sum\limits_{i>
j}y_i^*y_jy_i^*y_j=\sum\limits_{i< j}y_i^*y_jy_i^*y_j+
\sum\limits_{i< j}y_j^*y_iy_j^*y_i =\sum\limits_{i< j}
((y_i^*y_j)^2+ (y_j^*y_i)^2)\,.
\end{equation*}
Using the fact that $(y_i^*y_j)^2+(y_j^*y_i)^2\leq y_i^*y_jy_j^*y_i+
y_j^*y_iy_i^*y_j$\,, $1\leq i, j\leq d$\,, it follows that
\begin{eqnarray}\label{eq66776677}
\sum\limits_{i< j} ((y_i^*y_j)^2+ (y_j^*y_i)^2)&\leq & \sum\limits_{i< j} (y_i^*y_jy_j^*y_i+ y_j^*y_iy_i^*y_j)\\
&=& \sum\limits_{i< j} y_i^*y_jy_j^*y_i+ \sum\limits_{i> j} y_i^*y_jy_j^*y_i\nonumber\\
&=& \sum\limits_{i\ne j} y_i^*y_jy_j^*y_i\,.\nonumber
\end{eqnarray}
Therefore, we conclude that
\begin{eqnarray*}
(\text{Id}_n\otimes \mathbb{E})((Y^*Y)^2)&\leq & \left(\sum\limits_{i=1}^d y_i^*y_i\right)^2+ 2 \sum\limits_{i\ne j} y_i^*y_jy_j^*y_i\\
&\leq & \left(\sum\limits_{i=1}^d y_i^*y_i\right)^2 + 2 \sum\limits_{i, j=1}^d y_i^*y_jy_j^*y_i\\
&=& \left(\sum\limits_{i=1}^d y_i^*y_i\right)^2+ 2\sum\limits_{i=1}^d y_i^*\left(\sum\limits_{j=1}^d y_jy_j^*\right)y_i
\end{eqnarray*}
Recalling the definition (\ref{eq99900099900}), and using (\ref{eq11119}) we now obtain
\begin{eqnarray}\label{eq34343434}
(\text{Id}_n\otimes \mathbb{E})((Y^*Y)^2)&\leq & |||\{y_i\}_{i=1}^d|||^2\left(\sum\limits_{i=1}^d y_i^*y_i+ 2\sum\limits_{i=1}^d y_i^*y_i\right)\\
&=& 3|||\{y_i\}_{i=1}^d|||^2(\text{Id}_n\otimes \mathbb{E})(Y^*Y)\,.\nonumber
\end{eqnarray}
A similar proof based on (\ref{eq11120}) shows that
\begin{eqnarray}\label{eq34343435}
(\text{Id}_n\otimes \mathbb{E})((YY^*)^2)&\leq & 3|||\{y_i\}_{i=1}^d|||^2(\text{Id}_n\otimes \mathbb{E})(YY^*)\,.
\end{eqnarray}
Next we prove the following
\begin{lemma}\label{lem3333}
Let $x_1\,, \ldots \,,x_d\in M_n(\mathbb{C})$. Then there exists
$X\in M_n(L^\infty([0, 1]))$ such that \[ (\text{Id}_n\otimes
E)(X)=\sum\limits_{i=1}^d x_i\otimes e_i\,, \] satisfying, moreover,
\begin{eqnarray}\label{eq7777003}
\|X\|_{M_n(L^\infty([0, 1]))}&\leq & \sqrt{3}\max\left\{\left\|\sum\limits_{i=1}^d
x_i^*x_i\right\|^\frac12\,, \left\|\sum\limits_{i=1}^d
x_ix_i^*\right\|^\frac12\right\}\,.
\end{eqnarray}
\end{lemma}

As in the case of independent standard complex Gaussians, the
crucial point in the argument is the following version of Lemma
\ref{lem4567}, whose proof carries over verbatim to this setting,
except for choosing $C=\frac{\sqrt{3}}{2}$\,.

\begin{lemma}\label{lem12}
If $y_1\,, \ldots \,,y_d\in M_n(\mathbb{C})$ satisfy
\begin{eqnarray}\label{eq7777743}
\max\left\{\left\|\sum\limits_{i=1}^d
y_i^*y_i\right\|^\frac12\,, \left\|\sum\limits_{i=1}^d
y_iy_i^*\right\|^\frac12\right\}&=&1\,,
\end{eqnarray}
then there exists $Z\in M_n(L^\infty([0, 1]))$ such that
$\|Z\|_{M_n(L^\infty([0, 1]))}\leq \frac{\sqrt{3}}{2}$\,, and,
moreover, when $z_1\,, \ldots \,,z_d$ are defined by
$(\text{Id}_n\otimes E)(Z)=\sum\limits_{i=1}^d z_i\otimes e_i$\,,
then
\begin{eqnarray*}
\max\left\{\left\|\sum\limits_{i=1}^d
(y_i-z_i)^*(y_i-z_i)\right\|^\frac12\,, \left\|\sum\limits_{i=1}^d
(y_i-z_i)(y_i-z_i)^*\right\|^\frac12\right\}&\leq
&\frac12\,.
\end{eqnarray*}
\end{lemma}

Hence, for all $y\in M_n(\mathbb{C}^d)$ there is $Z\in
M_n(L^\infty([0, 1]))$ such that $\|Z\|_{M_n(L^\infty([0, 1]))}\leq
\frac{\sqrt{3}}{2}|||y|||_{M_n(H)}\,,$ and $|||(\text{Id}_n\otimes
E)(Z)-y|||_{M_n(\mathbb{C}^d)}\leq
 {\frac12}|||y|||_{M_n(\mathbb{C}^d)}$\,.
An application of Lemma \ref{lem25} with $V=M_n(L^\infty([0,
1]))$\,, $W=M_n(\mathbb{C}^d)$\,, $T=\text{Id}_n\otimes E$\,, the
map $\phi:W\rightarrow V$ be defined by $\phi(y)=Z$\,, $\forall y\in
W$\,, $C=\frac{\sqrt{3}}{2}$ and $\delta=\frac12$ shows that for all
$x\in M_n(\mathbb{C}^d)$ there exists $X\in M_n(L^\infty([0,1]))$
such that $(\text{Id}_n\otimes E)(X)=x$ and
\begin{equation}\label{eq5656567677}\|X\|_{M_n(L^\infty([0, 1]))}\leq
\sqrt{3}|||x|||_{M_n(\mathbb{C}^d)}\,. \end{equation} This completes
the proof of Lemma \ref{lem3333}\,. As explained before,
(\ref{eq5656567677}) implies that \[ \|(\text{Id}_n\otimes
q)(X)\|_{M_n({L^\infty([0,1])}/{\text{Ker}(E)})}\leq
\sqrt{3}\left|\left|\left|\{x_i\}_{i=1}^d\right|\right|\right|_{M_n(\mathbb{C}^d)}\,.
\] We conclude that there exists a linear bijection
$\widehat{E}:{L^\infty([0, 1])}/{\text{Ker}(E)}\rightarrow
\mathbb{C}^d$ such that $\widehat{E}(q(r_i))=e_i=E(r_i)$\,, for all
$1\leq i\leq d$\,, and moreover, with respect to the operator space
structure of the quotient space ${L^\infty([0,
1])}/{\text{Ker}(E)}$\,, the inequalities (\ref{eq2233445577}) hold.
This completes the proof of Proposition \ref{dualradem}.
\end{proof}

\begin{rem}\label{rem458978}\rm
Let $c_1\,, c_2$ denote the best constants in the inequalities
\begin{eqnarray}\label{eq800988905}
c_1|||\{x_i\}_{i=1}^d|||^*&\leq &\left\|\sum\limits_{i=1}^d
x_i\otimes r_i\right\|_{L^1([0, 1]; S_1^n)}\leq
c_2|||\{x_i\}_{i=1}^d|||^*\,,
\end{eqnarray}
where $d, n$ are positive integers, and $x_1\,, \ldots \,,x_d\in
M_n(\mathbb{C})$\,. Then the following estimates hold
\[ \frac1{\sqrt{3}}\leq c_1\leq \frac1{\sqrt{2}}\,, \quad c_2=1\,. \]
Indeed, the estimate $c_1\leq \frac1{\sqrt{2}}$ is a consequence of
Szarek's result (see \cite{Sz}) that the best constant in the
classical Khintchine inequalities for Rademachers is
$\frac1{\sqrt{2}}$\,, while the estimate $\frac1{\sqrt{3}}\leq c_1$
follows by Theorem \ref{th4}, which also shows that $c_2\leq 1$\,.
Since $\mathbb{E}(|r_1|)=1$\,, we deduce by taking $d=n=1$ and
$x_1=1$ in (\ref{eq800988905}) that $c_2\geq 1$\,. Hence $c_2=1$\,.
\end{rem}

\section{A noncommutative Khintchine-type inequality for subspaces of $R\oplus C$} \setcounter{equation}{0}
Let $H\subseteq R\oplus C$ be a subspace, equipped with the Hilbert
space structure induced by the usual direct sum of Hilbert spaces
inner product. More precisely, given $\xi\in H$\,, write
$\xi=({\xi}_R, {\xi}_C)\in R\oplus C$\,;  then
\begin{eqnarray}\label{eq1}
{\langle \xi, \eta\rangle}_H &=& {\langle {\xi}_R,
{\eta}_R\rangle}_R+ {\langle {\xi}_C, {\eta}_C\rangle}_C\,, \quad
\forall \xi, \eta\in H\,.
\end{eqnarray}
Consider $R\oplus C$ equipped with the operator space structure of
the $l_\infty$-direct sum $R{\oplus}_\infty C$\,. Note that the norm
induced on $H$ by the inner product $\langle \cdot, \cdot \rangle_H$
is not the same as the one coming from $R{\oplus}_\infty C$\,. For
all $\xi\in H$\,, define further
\[ U_1(\xi)={\xi}_R\,, \quad U_2(\xi)={\xi}_C\,. \]
Then $U_1\in {\mathcal B}(H, R)$\,, respectively $U_2\in {\mathcal
B}(H, C)$ and formula (\ref{eq1}) becomes
\begin{eqnarray}\label{eq2}
{\langle \xi, \eta\rangle}_H &=& {\langle U_1({\xi}),
U_1({\eta})\rangle}_R+ {\langle U_2({\xi}),
U_2({\eta})\rangle}_C\,, \quad \forall \xi, \eta\in H.
\end{eqnarray}
The operator $U: H\rightarrow R\oplus C$ defined by $U=\left( \begin{array} [c] {c}%
U_1\\U_2
\end{array} \right)$
is an isometry, where $H$ and $R\oplus C$ are equipped with the
above Hilbert space structure. This implies that
$U_1^*U_1+U_2^*U_2=I$\,, where $I$ denotes the identity operator on
$H$\,. Let
\begin{eqnarray}\label{eq111}
A&=&U_2^*U_2\in {\mathcal B}(H)\,. \end{eqnarray} Then $0\leq A\leq
I$\,.

We now discuss the operator space structure of $H$. Let $n$ be a
positive integer. Then for all $r\in \mathbb{N}$\,, all $x_i\in
M_n(\mathbb{C})$ and all ${\xi}_i\in H$\,, $1\leq i\leq r$\,, we
have
\begin{eqnarray}\label{eq12}
\left\|\sum\limits_{i=1}^r x_i\otimes {\xi}_i
\right\|_{M_n({H})}&=&\max\left\{ \left\|\sum\limits_{i=1}^r
x_i\otimes U_1{\xi}_i\right\|_{M_n(R)}\,, \left\|\sum\limits_{i=1}^r
x_i\otimes U_2{\xi}_i\right\|_{M_n(C)}\right\}\,.
\end{eqnarray}
We claim that
\begin{eqnarray}\label{eq3}
\left\|\sum\limits_{i=1}^r x_i\otimes {\xi}_i
\right\|_{M_n({H})}&=&\max\left\{ \left\|\sum\limits_{i, j=1}^r
\langle (I-A){\xi}_i, {\xi}_j\rangle_H x_ix_j^*\right\|^\frac12\,,
\left\|\sum\limits_{i, j=1}^r \langle A{\xi}_i, {\xi}_j\rangle_H
x_i^*x_j\right\|^\frac12 \right\}\,.
\end{eqnarray}
Indeed, by the definition of operator space matrix norms on $R$ and $C$ we have
\begin{equation*}
\left\|\sum\limits_{i=1}^r x_i\otimes U_1{\xi}_i\right\|_{M_n(R)}=
\left\|\sum\limits_{i, j=1}^r x_ix_j^*\langle U_1{\xi}_i,
U_1{\xi}_j\rangle_R\right\|^\frac12=\left\|\sum\limits_{i, j=1}^r
x_ix_j^*\langle (I-A){\xi}_i, {\xi}_j\rangle_H\right\|^\frac12\,,
\end{equation*}
respectively,
\begin{equation*}
\left\|\sum\limits_{i=1}^r x_i\otimes
U_2{\xi}_i\right\|_{M_n(C)}=\left\|\sum\limits_{i, j=1}^r
x_i^*x_j\langle U_2{\xi}_i, U_2{\xi}_j\rangle_C\right\|^\frac12
=\left\|\sum\limits_{i, j=1}^r x_i^*x_j\langle A{\xi}_i,
{\xi}_j\rangle_H\right\|^\frac12\,,
\end{equation*}
and the claim is proved.

Let ${\mathcal A}$ be the CAR algebra over the Hilbert space $H$.
Recall that ${\mathcal A}$ is a unital $C^*$-algebra (unique up to $*$-isomorphism) with the property
that there exists a linear map
\[ H\ni f\mapsto a(f)\in {\mathcal A} \]
whose range generates ${\mathcal A}\,,$ satisfying for all $f,
g\in H$ the anticommutation relations
\begin{eqnarray}\label{eq4}
a(f)a(g)^*+a(g)^*a(f)&=&\langle f, g\rangle_H I\\
a(f)a(g)+a(g)a(f)&=&0\,.\nonumber
\end{eqnarray}
Let $\omega_A$ be the gauge-invariant quasi-free state on ${\mathcal
A}$ corresponding to the operator $A$ ($0\leq A\leq I$) associated
to the subspace $H$ of $R\oplus C$. Recall that a state $\omega$ on
${\mathcal A}$ is called gauge-invariant if it is invariant under
the group of gauge transformations $\tau_\theta
(a(f))=a(e^{i\theta}f)$\,, $\forall \theta\in [0, 2\pi)$\,. It turns
out (see  \cite{BV} and \cite{BR}) that a gauge-invariant quasi-free
state $\omega$ on ${\mathcal A}$ is completely determined by one
truncated function $\omega_T$. More precisely, a functional
$\omega_T(\cdot, \cdot)$ over the monomials in $a^*(f)$ and $a(g)\,,
\forall f, g\in H$\,, which is linear in the first argument and
conjugate-linear in the second determines a gauge-invariant
quasi-free state $\omega$ on ${\mathcal A}$ if and only if
\begin{equation}\label{eq5}
0\leq \omega_T(a(f)^*, a(f))\leq \|f\|^2\,, \quad\forall f\in H\,.
\end{equation}
Now, given the operator $0\leq A\leq I$, define
\[  \omega_T^{A}(a(f)^*, a(g)):=\langle Ag, f\rangle_H\,. \]
The positivity condition (\ref{eq5}) is clearly satisfied. Let
$\omega_A$ be the gauge-invariant quasi-free state on ${\mathcal
A}$ determined by the truncated function $\omega_T^A$. Then for
all $n\geq 1$, the $n$-point functions of $\omega_A$ have the form
\begin{equation}\label{eq6}
\omega_A(a(f_n)^*\ldots a(f_1)^*a(g_1)\ldots a(g_m))=\delta_{nm}
\text{det}(\langle Ag_i, f_j\rangle_H\,, i, j)\,, \quad\forall f_1,
\ldots \,,f_n, g_1, \ldots \,,g_m\in H\,.
\end{equation}
Given $b\in {\mathcal A}$, the map
\[ H\ni f\mapsto \omega_A(a(f)b^*+b^*a(f))\in \mathbb{C} \]
is a bounded linear functional on $\A$. By the Riesz representation
theorem, there exists a unique element $E_A(b)\in H$ such that
\begin{eqnarray}\label{eq7}
\langle f, E_A(b)\rangle_H&=&\omega_A(a(f)b^*+b^*a(f))\,,
\quad\forall f\in H\,.
\end{eqnarray}
Equivalently,
\begin{eqnarray}\label{eq8}
\langle E_A(b), f\rangle_H&=&\omega_A(b{a(f)}^*+{a(f)}^*b)\,,
\quad\forall f\in H\,.
\end{eqnarray}
We obtain in this way a bounded linear map $E_A:\A\rightarrow H$\,.
By uniqueness in the Riesz representation theorem and the
anticommutation relations (\ref{eq4}) it follows that
\begin{eqnarray}\label{eq9}
E_A(a(f))=f\,, \quad \forall f\in H\,.
\end{eqnarray}

Consider the GNS representation $(\pi_{\omega_A}, H,
\xi_{\omega_A})$ associated to $(\A, \omega_A)$. For simplicity of
notation, write $\pi_{\omega_A}=\pi_A$ and $\xi_{\omega_A}=\xi_A$
(the cyclic unit vector for the representation). Then for all $f\in
H$ and all $b\in \A$\,,
\[ \omega_A(a(f)b^*+b^*a(f))=\langle \pi_A(a(f)b^*+b^*a(f))\xi_A\,,
\xi_A\rangle_H=\langle \{\pi_A(a(f))\,, \pi_A(b^*)\}\xi_A\,,
\xi_A\rangle_H\,, \] where $\{K, L\}=KL+LK$\,. Equivalently,
\begin{eqnarray*}
\omega_A(b{a(f)}^*+{a(f)}^*b)&=&\langle \{\pi_A(a(f)^*)\,,
\pi_A(b)\}\xi_A\,, \xi_A\rangle_H\,, \quad \forall f\in H\,,
\forall b\in \A\,.
\end{eqnarray*}
Note that the map
\[ \A\ni c\mapsto \langle\{\pi_A(a(f)^*)\,, c\}\xi_A\,,
\xi_A\rangle_H\in \mathbb{C} \] extends to a normal (positive)
linear functional on the von Neumann algebra
$\overline{\pi_A(\A)}^{\text{sot}}$. This implies that $E_A$ extends
to a bounded linear map on the von Neumann algebra generated by
$\pi_A(\A)$ and moreover the range of the dual map $E_A^*$ is
contained in the predual of $\overline{\pi_A(\A)}^{\text{sot}}$.\\

%Let $r\in \mathbb{N}$. Denote by $\{f_i\}_{1\leq r\leq r}$ the canonical unit basis in ${\mathbb{C}}^r$. Let
%\[ a_i:=a(f_i)\in {\A}\,, \quad \forall 1\leq i\leq r. \] By (\ref{eq8}) and the anticommutation relations (\ref{eq4}) it follows that
%\[ E(a_i)=f_i\,, \quad \forall 1\leq i\leq r\,. \]
%Furthermore, define
%\[ A(i, j):=\langle Af_i, f_j\rangle_H\,, \quad \forall 1\leq i, j\leq r\,. \]

With the notation set forth above, we prove the following

\begin{theorem}\label{th1}
The map $E_A:{\A}\rightarrow H$ yields a complete isomorphism
\[ H\cong {{\A}/{\text{Ker}(E_A)}} \]
with cb-isomorphism constant $\leq \sqrt{2}$\,. More precisely, if
$q_A:{\A}\rightarrow {{\A}/{\text{Ker}(E_A)}}$ denotes the quotient
map, then given any positive integers $n, r$ we have for all $x_i\in
M_n(\mathbb{C})$ and $b_i\in \A$\,, $1\leq i\leq r$:
\begin{equation}\label{eq13}
\left\|\sum\limits_{i=1}^{r} x_i\otimes
E_A(b_i)\right\|_{M_n(H)}\leq \left\|\sum\limits_{i=1}^r x_i\otimes
q_A(b_i)\right\|_{M_n({{\A}/{\text{Ker}(E_A)}})}\leq
\sqrt{2}\left\|\sum\limits_{i=1}^r x_i\otimes
E_A(b_i)\right\|_{M_n(H)}\,.
\end{equation}
Furthermore, the dual map $E_A^*$ is a complete isomorphism of $H^*$
onto a subspace of the predual of
$\overline{\pi_A(\A)}^{\text{sot}}$.
\end{theorem}

\begin{rem}\label{re667788}\rm Note that Theorem \ref{th1} is
equivalent to the statement that for any positive integers $n$\,,
$r$ we have for all $x_i\in M_n(\mathbb{C})$ and $\xi_i\in H$\,,
$1\leq i\leq r$\,,
\begin{equation}\label{eq131313}
\left\|\sum\limits_{i=1}^{r} x_i\otimes {\xi}_i\right\|_{M_n(H)}\leq
\left\|\sum\limits_{i=1}^r x_i\otimes
q_A(a(\xi_i))\right\|_{M_n({{\A}/{\text{Ker}(E_A)}})}\leq
\sqrt{2}\left\|\sum\limits_{i=1}^r x_i\otimes
\xi_i\right\|_{M_n(H)}\,.
\end{equation}
\end{rem}
Indeed, to prove that (\ref{eq13}) implies (\ref{eq131313}), put
$b_i:=a(\xi_i)$\,, $1\leq i\leq r$ and use the fact that by
(\ref{eq9}), $E_A(a(\xi_i))=\xi_i$\,, $1\leq i\leq r$\,. To prove
that, conversely, (\ref{eq131313}) implies (\ref{eq13}), put
$\xi_i:=E_A(b_i)$\,, $1\leq i\leq r$\,. Then
$E_A(b_i-a(\xi_i))=0$\,, which implies that $q_A(b_i-a(\xi_i))=0$\,,
so the middle term of (\ref{eq13}) is equal to the middle term of
(\ref{eq131313})\,. The equivalence of (\ref{eq13}) and
(\ref{eq131313}) will be used several times in the following.\\

\noindent {\em Proof of Theorem \ref{th1}}. We first prove the
theorem in the finite dimensional case.

Assume $\text{dim}(H)=d< \infty$. Consider the associated operator
$A$ ($0\leq A\leq I$) defined by (\ref{eq111}). There exists an
orthonormal basis $\{e_i\}_{1\leq i\leq d}$ of $H$  with respect to
which the matrix $A$ is diagonal. That is,
\begin{equation}\label{eq77777777789778}
\langle A e_i\,, e_j\rangle_H=\nu_i\delta_{ij}\,, \quad \forall
1\leq i, j\leq d\,, \end{equation} which implies that $0\leq
\nu_i\leq 1$\,, $\forall 1\leq i\leq d$\,.

\noindent Let $\A$ be the CAR-algebra over $H$ and $\omega_A$ be the
quasi-free state on $\A$ corresponding to the operator $A$\,.
Further, set
\begin{eqnarray*} a_i&:=& a(e_i)\,, \quad \forall 1\leq i\leq d\,. \end{eqnarray*}
By (\ref{eq6}) it follows that
\begin{eqnarray}\label{eq56677889901}
\omega_A(a_i^*a_j)&=&\nu_i\delta_{ij}\,, \quad \forall 1\leq i,
j\leq d\,,
\end{eqnarray}
and, respectively,
\begin{eqnarray}\label{eq56677889902}
\omega_A(a_ia_j^*)&=&(1-\nu_i)\delta_{ij}\,, \quad \forall 1\leq i, j\leq d\,.
\end{eqnarray}

Let $n$ be a positive integer. Given $x_1\,, \ldots \,,x_d\in
M_n(\mathbb{C})$\,, we have by (\ref{eq3}) that
\begin{equation}\label{eq223322222}
\left|\left|\left|\{x_i\}_{i=1}^d\right|\right|\right|_A:=\left\|\sum\limits_{i=1}^{d}
x_i\otimes {e_i}\right\|_{M_n(H)}=\max\left\{
\left\|\sum\limits_{i=1}^{d} (1-\nu_i)x_ix_i^*\right\|^\frac12\,,
\left\|\sum\limits_{i=1}^{d} \nu_i x_i^*x_i\right\|^\frac12
\right\}\,. \end{equation} In view of Remark \ref{re667788}, we have
to prove that
\begin{eqnarray}\label{eq10}
\qquad \quad\max\left\{ \left\|\sum\limits_{i=1}^{d}
(1-\nu_i)x_ix_i^*\right\|^\frac12\,, \left\|\sum\limits_{i=1}^{d}
\nu_i
x_i^*x_i\right\|^\frac12 \right\}&\!\!\leq \!\!&\left\|\sum\limits_{i=1}^{d} x_i\otimes q_A(a_i)\right\|_{M_n({{\A}/{\text{Ker}(E_A)}})}\\
&\!\!\leq \!\!& \sqrt{2}\max\left\{ \left\|\sum\limits_{i=1}^{d}
(1-\nu_i)x_ix_i^*\right\|^\frac12\,, \left\|\sum\limits_{i=1}^{d}
\nu_i x_i^*x_i\right\|^\frac12 \right\}.\nonumber
\end{eqnarray}
We first prove the left hand side inequality in (\ref{eq10})\,. For each $1\leq i\leq d$ set
\begin{eqnarray}\label{eq88}
\phi_i^A(b)&:=& \omega_A(a_i^*b+ba_i^*)\,, \quad \forall b\in \A\,.
\end{eqnarray}
Note that $\langle E_A(b), e_i\rangle_H=\phi_i^A(b)$\,, for all $b\in\A$ and that by the anticommutation relations (\ref{eq4})\,,
\[ \phi_i^A(a_j)=\delta_{ij}\,, \quad \forall 1\leq i, j\leq d\,. \]
In particular, $E_A(a_i)=e_i$\,, $\forall 1\leq i\leq d$\,.

\begin{lemma}\label{lem677767776} For all $1\leq i\leq d$ we have
\begin{eqnarray}\label{eq38}
\omega_A(a_i^*b)={\nu_i}\phi_i^A(b)\,, &&
\omega_A(ba_i^*)=(1-\nu_i)\phi_i^A(b)\,, \quad\forall b\in \A\,.
\end{eqnarray}
\end{lemma}

\begin{proof}
We consider a special representation of the CAR algebra $\A$\,. Let
$e=\left(\begin{array}
[c]{cc}%
0 & 1\\
0 & 0
\end{array}
\right)$, $u=\left(\begin{array}
[c]{cc}%
1 & 0\\
0 & -1
\end{array}
\right),$ $I_2=I_{M_2(\mathbb{C})}$ and set
\begin{equation}\label{eq300004509}
a_1^{\prime}:=e\otimes (\otimes_{j=2}^d I_2)\,, \quad
a_i^{\prime}:=\left(\otimes_{j=1}^{i-1} u\right)\otimes e\otimes
(\otimes_{j=i+1}^d I_2)\,, \quad 2\leq i\leq d\,. \end{equation}
Since $u^2=I_2$\,, $ee^*+e^*e=I_2$\,, $eu+ue=0$\,, it follows that
$\{a_i^{\prime}\}_{1\leq i\leq d}$ satisfy the CAR relations
(\ref{eq4}). Thus $C^*(\{a_1^{\prime}\,, \ldots \,, a_d^{\prime}\})=
\otimes_{i=1}^d M_2(\mathbb{C})$ (see \cite{BR} and \cite{D}), and
there is a $*$-isomorphism $\psi: \A\rightarrow
C^*(\{a_1^{\prime}\,, \ldots \,, a_d^{\prime}\})$ such that
$\psi(a_i)=a_i^\prime$\,, $\forall 1\leq i\leq d$\,. From now on we
identify $\A$ with $\otimes_{i=1}^d M_2(\mathbb{C})$\,, and write
$a_i^{\prime}=a_i$\,, $1\leq i\leq d$\,. Then, by \cite{PS} (see pp.
4 and 5),
\begin{equation}\label{eq40004004}
\omega_{A}(b):= \left(\otimes_{i=1}^d \psi_i\right)(b)\,, \quad
\forall b \in \A\,,
\end{equation}
where $\psi_i(h)= \text{Tr}\left(\left(\begin{array}
[c]{cc}%
1-\nu_i & 0\\
0 & \nu_i
\end{array}
\right)h\right)$\,, $\forall h\in M_2(\mathbb{C})$\,, $1\leq i\leq
d$\,.

We first show that for all $1\leq i\leq d$\,,
\begin{eqnarray}\label{eq912}
(1-\nu_i)\omega_A((a_i)^*b)&=&{\nu_i}\omega_A(b(a_i)^*)\,, \quad
\forall b\in \A\,.
\end{eqnarray}
To check (\ref{eq912}), it is enough to look at simple tensors
$b=b_1\otimes b_2\otimes \ldots \otimes b_d \in \A$\,. Consider
first the case $i=1$\,. Then
\[ \omega_{A}
((a_1)^*b)=\psi_1(e^*b_1)\prod\limits_{i=2}^d \psi_i(b_i)\,, \quad
\omega_{A}(b(a_1)^*)=\psi_1(e^*b_1)\prod_{i=2}^d \psi_i(b_i)\,.
\] Let $b_1=\left(\begin{array}
[c]{cc}%
b_1^{(11)} & b_1^{(12)}\\
b_1^{(21)} & b_1^{(22)}
\end{array}
\right)$\,. Then $\psi_1(e^*b_1)=\psi_1\left(\left(\begin{array}
[c]{cc}%
0 & 0\\
b_1^{(11)} & b_1^{(12)}
\end{array}
\right)\right)=\text{Tr}\left(\left(\begin{array}
[c]{cc}%
0 & 0\\
0 & \nu_1 b_1^{(12)}
\end{array}
\right)\right)={\nu_1}b_1^{(12)}$, respectively,
$\psi_1(b_1e^*)=\pi_1\left(\left(\begin{array}
[c]{cc}%
b_1^{(12)} & 0\\
b_1^{(22)} & 0
\end{array}
\right)\right)=\text{Tr}\left(\left(\begin{array}
[c]{cc}%
(1-\nu_1)b_1^{(12)} & 0\\
{\nu_1}b_1^{(22)} & 0
\end{array}
\right)\right)=(1-\nu_1)b_1^{(12)}$\,. Hence
\[ \omega_A((a_1)^*b)={\nu_1}b_1^{(12)}\prod\limits_{i=2}^d
\psi_i(b_i)\,, \quad
\omega_A(b(a_1)^*)=(1-\nu_1)b_1^{(12)}\prod\limits_{i=2}^d
\psi_i(b_i)\,, \] which imply (\ref{eq912}). The case when $i\neq 1$
can be proved in a similar way, using the fact that for all
$b=b_1\otimes b_2\otimes \ldots \otimes b_d \in \A$\,, we have
$ub_j=b_ju$\,, $\forall 1\leq j\leq i-1$\,.

Then, for $1\leq i\leq d$ we deduce by (\ref{eq912}) that for all
$b\in \A$\,, we have
\begin{eqnarray*}
{\nu_i}\omega_A((a_i)^*b+ b(a_i)^*)=
{\nu_i}\omega_A((a_i)^*b)+(1-\nu_i)\omega_A((a_i)^*b)
=\omega_A((a_i)^*b)\,,
\end{eqnarray*}
and, respectively,
\begin{eqnarray*}(1-\nu_i)\omega_A((a_i)^*b+ b(a_i)^*)
= {\nu_i}\omega_A(b(a_i)^*)+(1-\nu_i)\omega_A(b(a_i)^*)
=\omega_A(b(a_i)^*)\,.\end{eqnarray*} The proof is complete.
\end{proof}

\begin{lemma}\label{lem4141}
Let $X\in M_n(\A)$\,. By letting
\begin{eqnarray}\label{eq233}
x_i&:=& (\text{Id}_n\otimes \phi_i^A)(X)\,, \quad \forall 1\leq i\leq d\,,
\end{eqnarray}
we have
\begin{eqnarray}\label{eq25}
(\text{Id}_n\otimes E_A)(X)&=& \sum\limits_{i=1}^d x_i\otimes e_i\,.
\end{eqnarray}
Then, with the above notation it follows that
\begin{eqnarray}\label{eq36}
\|X\|_{M_n({\A})}&\geq &\max\left\{\left\|\sum\limits_{i=1}^{d}
\nu_i x_i^*x_i\right\|^{\frac12}\,, \left\|\sum\limits_{i=1}^{d}
(1-\nu_i)x_ix_i^*\right\|^\frac12\right\}\,.
\end{eqnarray}
\end{lemma}

\begin{proof}
Let $X\in M_n(\A)$\,. Then $X$ is of the form $X=\sum\limits_{j=1}^r
y_j\otimes b_j$\,, where $r\in \mathbb{N}$\,, $y_j\in
M_n(\mathbb{C})$ and $b_j\in \A$\,, $1\leq j\leq r$\,. For all
$1\leq i\leq d$\,, let $x_i$ be defined by (\ref{eq233})\,. Then
\begin{equation}\label{eq776766779} x_i=\sum\limits_{j=1}^r
\phi_i^A(b_j)y_j\,, \quad 1\leq i\leq d\,.
\end{equation}
Further set
\begin{equation}\label{eq7767667799}
Z:=\sum\limits_{i=1}^d x_i\otimes a_i\in M_n(\A)\,.
\end{equation}
To each state $\omega$ on $M_n(\mathbb{C})$ we can associate a
positive sesquilinear form on $M_n(\A)$ given by
\[ s_\omega(c, d):=(\omega\otimes \omega_A)(d^*c)\,, \quad \forall
c, d\in M_n(\A)\,. \] By (\ref{eq776766779}), (\ref{eq7767667799})
and (\ref{eq38}), we obtain
\begin{eqnarray*}
s_\omega(X, Z)=\sum\limits_{i=1}^d \sum\limits_{j=1}^r
\omega(x_i^*y_j)\omega_A(a_i^*b_j)&=&\sum\limits_{i=1}^d
\sum\limits_{j=1}^r \nu_i \omega(x_i^*y_j)\phi_i^A(b_j)\\&=&
\sum\limits_{i=1}^d \nu_i \omega(x_i^*x_i)=s_\omega(Z, Z)\,,
\end{eqnarray*}
where the last equality follows from (\ref{eq56677889901})\,. Hence
$s_\omega(X-Z, Z)=0$\,, and therefore
\[ s_\omega(X, X)=s_\omega(Z, Z)+s_\omega(X-Z, X-Z)\geq s_\omega(Z,
Z)\,. \] It follows that
\[ \omega\left(\sum\limits_{i=1}^d \nu_ix_i^*x_i\right)=s_\omega(Z,
Z)\leq s_\omega (X, X)\leq \|X\|^2\,, \] for every state $\omega$ on
$M_n(\mathbb{C})$\,, and hence
\[ \left\|\sum\limits_{i=1}^d \nu_ix_i^*x_i\right\|\leq \|X\|^2\,.
\]
The same argument applied to the positive sesquilinear form
\[ s_\omega^{\prime}(c, d):=(\omega\otimes \omega_A)(cd^*)\,, \quad
\forall c, d\in M_n(\A) \] gives by (\ref{eq56677889902}) that
\[ \omega\left(\sum\limits_{i=1}^d (1-\nu_i)x_ix_i^*\right)=s_\omega^{\prime}(Z,
Z)\leq s_\omega^{\prime} (X, X)\leq \|X\|^2\,, \] for every state
$\omega$ on $M_n(\mathbb{C})$\,, and hence
\[ \left\|\sum\limits_{i=1}^d (1-\nu_i)x_ix_i^*\right\|\leq
\|X\|^2\,.
\]
This completes the proof.
\end{proof}

\begin{rem}\label{rem5}\rm
For all $X\in M_n(\A)$ we have
\begin{eqnarray}\label{eq77777}
\|(\text{Id}_n\otimes E_A)(X)\|_{M_n(H)}&\leq &
\|(\text{Id}_n\otimes q_A)(X)\|_{M_n({{\A}/{\text{Ker}(E_A)})}}\,.
\end{eqnarray}
This follows by a similar argument as the one used to prove
(\ref{eq22334460}). In particular, given $x_1, \ldots , x_d\in
M_n(\mathbb{C})$\,, by letting $X=\sum\limits_{i=1}^d x_i\otimes
a_i\in M_n(\A)$\,, an application of (\ref{eq77777}) yields the left
hand side inequality in (\ref{eq10})\,.
\end{rem}

We now prove the right hand side inequality in (\ref{eq10})\,. Let $y_1\,, \ldots \,, y_{d}\in M_n(\mathbb{C})$\,. Set
\begin{eqnarray*} Y&:=&\sum\limits_{i=1}^{d}y_i\otimes a_i\in M_n(\A)\,. \end{eqnarray*}
We will compute $(\text{Id}_n\otimes \omega_{A})(Y^*Y)$\,,
$(\text{Id}_n\otimes \omega_{A})(YY^*)$\,, $(\text{Id}_n\otimes
\omega_{A})((Y^*Y)^2)$ and
$(\text{Id}_n\otimes \omega_{A})((YY^*)^2)$\,. We have
\[ Y^*Y=\sum\limits_{i, j=1}^{d} y_i^*y_j\otimes a_i^*a_j\,, \quad YY^*=\sum\limits_{i, j=1}^{d} y_iy_j^*\otimes a_ia_j^*\,. \] By (\ref{eq56677889901}) and (\ref{eq56677889902})
it follows immediately that
\begin{eqnarray}
(\text{Id}_n\otimes \omega_{A})(Y^*Y)&=& \sum\limits_{i=1}^{d} \nu_iy_i^*y_i\label{eq19}\\
(\text{Id}_n\otimes \omega_{A})(YY^*)&=& \sum\limits_{i=1}^{d}
(1-\nu_i)y_iy_i^*\,.\label{eq20}
\end{eqnarray}
Furthermore, in order to compute $(\text{Id}_n\otimes
\omega_{A})((Y^*Y)^2)$\,, note that
\begin{eqnarray}\label{eq2222222222}
Y^*Y&=& \sum\limits_{i, j=1}^d y^*_iy_j\otimes (a^*_ia_j-\delta_{ij}\nu_iI)+\sum\limits_{i=1}^d \nu_iy_i^*y_i\otimes I\,.
\end{eqnarray}
Consider the vectors
\begin{eqnarray*}
f_{ij}&:=& a_i^*a_j-\delta_{ij}\nu_iI\,, \quad \forall 1\leq i, j\leq d\,.
\end{eqnarray*}
We claim that $\{I, f_{ij}\,, 1\leq i, j\leq d\}$ is an orthogonal
set in $L^2(\A)$ with respect to the positive sesquilinear form on $\A$ given by
${\A\times \A}\ni (c, d)\mapsto \omega_A(d^*c)\in \mathbb{C}$\,,
satisfying $\omega_A(I)=1$ and
\begin{eqnarray}\label{eq545454546}
\omega_A(f_{ij}^*f_{ij})&=&\nu_j(1-\nu_i)\,, \quad \forall 1\leq i, j\leq d\,. \end{eqnarray}
Indeed, for $1\leq i, j\leq d$\,,
\begin{eqnarray*}
\omega_A(f_{ij}^*f_{ij})&=&\omega_A(a_j^*a_ia_i^*a_j)-{\nu_i}\omega_A(a_j^*a_i+a_i^*a_j)\delta_{ij}+\nu_i^2\delta_{ij}\,.
\end{eqnarray*}
By the anticommutation relations (\ref{eq4}), together with (\ref{eq6}) we get
\[
\omega_A(a_j^*a_ia_i^*a_j)=\omega_A(a_j^*(I-a_i^*a_i)a_j)
=\omega_A(a_j^*a_j)-\omega_A(a_j^*a_i^*a_ia_j) =
\nu_j-\nu_i\nu_j(1-\delta_{ij})=\nu_j(1-\nu_i)+\nu_i\nu_j\delta_{ij}\,,
\] wherein we have also used the fact that $a_i^2=0$\,, $1\leq i\leq
d$\,. Furthermore,
$\omega_A(a_j^*a_i+a_i^*a_j)=2\nu_i\delta_{ij}$\,. Hence
$\omega_A(f_{ij}^*f_{ij})=
\nu_j(1-\nu_i)+(\nu_i\nu_j-\nu_i^2)\delta_{ij} =\nu_j(1-\nu_i)$\,,
so (\ref{eq545454546}) is proved.

We now prove the orthogonality property of the set of vectors $\{I, f_{ij}\,, 1\leq i, j\leq d\}$\,.\\
First, note that for $1\leq i, j\leq d\,,$
\begin{eqnarray}\label{eq65656567}
\omega_A(f_{ij})&=& \omega_A(a_i^*a_j)-\nu_i\delta_{ij}=\nu_i\delta_{ij}-\nu_i\delta_{ij}=0\,.
\end{eqnarray}
It remains to show that for $1\leq i, j, k, l\leq d\,,$
\begin{eqnarray}\label{eq65656568}
\omega_A(f_{ij}^*f_{kl})=0\,, \quad \text{whenever} \,\,\,(i, j)\neq
(k, l)\,.
\end{eqnarray}
We have $f_{ij}^*f_{kl}=
a_j^*a_ia_k^*a_l-\nu_ka_j^*a_i\delta_{kl}-\nu_ia_k^*a_l\delta_{ij}+\nu_i\nu_k\delta_{ij}\delta_{kl}$\,.
We distinguish the following cases:
\[ 1) \,i=j\neq k=l\,, \,\,\,\,2) \,i\neq j\,, k=l\,, \,\,\,\,3) \,i=j\,, k\neq l\,, \,\,\,\,4) \,i\neq j\,, k\neq l\,, (i, j)\neq (k, l)\,. \]
Assume $ 1) \,i=j\neq k=l$. Then
\begin{eqnarray*}
\omega_A(f_{ii}^*f_{kk})&=& \omega_A(a_i^*a_ia_k^*a_k)-\nu_k\omega_A(a_i^*a_i)-\nu_i\omega_A(a_k^*a_k)+\nu_i\nu_k\\
&=&\omega_A(a_i^*(-a_k^*a_i)a_k)-\nu_k\nu_i-\nu_i\nu_k+\nu_i\nu_k\\
&=&-\omega_A(a_i^*a_k^*(-a_ka_i))-\nu_i\nu_k\\
&=&\nu_i\nu_k-\nu_i\nu_k=0
\end{eqnarray*}
Cases $2)$ and $3)$ are similar, so we only prove one of them. Assume $2) \,i\neq j\,, k=l$\,. Then
\[
\omega_A(f_{ij}^*f_{kk})=\omega_A(a_j^*a_ia_k^*a_k)-\nu_k\omega_A(a_j^*a_i)
=\omega_A(a_j^*(I\delta_{ik}-a_k^*a_i)a_k)=\omega_A(a_j^*a_k\delta_{ik})-\omega_A(a_j^*a_k^*a_ia_k\delta_{ik})
=0\,.
\] Respectively, assume $4) \,i\neq j\,, k\neq l\,, (i, j)\neq (k,
l)$\,. In this case,
$\omega_A(f_{ij}^*f_{kl})=\omega_A(a_j^*a_ia_k^*a_l)$\,. By
considering further the two possible subcases $4a)\, i\neq k$ and
$4b)\, i=k, j\neq l$\,, we deduce by (\ref{eq4}) and (\ref{eq6})
that $\omega_A(a_j^*a_ia_k^*a_l)=0$\,.

Then, based on the expansion (\ref{eq2222222222}) of $Y^*Y$ in terms
of the vectors $\{I, f_{ij}\,, 1\leq i, j\leq d\}$\,, we now get
\begin{eqnarray}
(\text{Id}_n\otimes \omega_A)((Y^*Y)^2)
&=&\sum\limits_{i, j=1}^d \nu_j(1-\nu_i)(y_i^*y_j)^*y_i^*y_j+\left(\sum\limits_{i=1}^d \nu_i y_i^*y_i\right)^2\label{eq19191919}\\
&=&\sum\limits_{j=1}^d \nu_jy_j^*\left(\sum\limits_{i=1}^d (1-\nu_i)y_iy_i^*\right)y_j+ \left(\sum\limits_{i=1}^d \nu_i y_i^*y_i\right)^2\nonumber\\
&\leq & \left(\sum\limits_{i=1}^d \nu_iy_i^*y_i\right)\left(\left\|\sum\limits_{i=1}^d (1-\nu_i)y_iy_i^*\right\|+\left\|\sum\limits_{i=1}^d \nu_iy_i^*y_i\right\|\right)\nonumber\\
&=&\left(\left\|\sum\limits_{i=1}^d \nu_iy_i^*y_i\right\|+\left\|\sum\limits_{i=1}^d (1-\nu_i)y_iy_i^*\right\|\right)(\text{Id}_n\otimes \omega_A)(Y^*Y)\,,\nonumber
\end{eqnarray}
where the last equality is given by (\ref{eq19})\,.

In order to estimate the term $(\text{Id}_n\otimes \omega_A)((YY^*)^2)$\,, note that
\begin{eqnarray}\label{eq2222222223}
YY^*&=&\sum\limits_{i, j=1}^d y_iy_j^*\otimes (a_ia_j^*-\delta_{ij}(1-\nu_i)I)+ \sum\limits_{i=1}^d (1-\nu_i)y_iy_i^*\otimes I\,.
\end{eqnarray}
We now consider the vectors
\begin{eqnarray*}
g_{ij}&:=& a_ia_j^*-\delta_{ij}(1-\nu_i)I\,, \quad \forall 1\leq i, j\leq d\,.
\end{eqnarray*}
With a similar proof it can be shown that $\{I, g_{ij}\,, 1\leq i, j\leq d\}$
is an orthogonal
set in $L^2(\A)$ with respect to the positive sesquilinear form on $\A$ given by
${\A\times \A}\ni (c, d)\mapsto \omega_A(cd^*)\in \mathbb{C}$\,,
satisfying
\begin{eqnarray}\label{eq54545454689}
\omega_A(g_{ij}g_{ij}^*)&=&\nu_j(1-\nu_i)\,, \quad \forall 1\leq i,
j\leq d\,. \end{eqnarray} Thus, based on the expansion
(\ref{eq2222222223}) of $YY^*$ in terms of the vectors $\{I,
g_{ij}\,, 1\leq i, j\leq d\}$\,, we obtain
\begin{eqnarray}
(\text{Id}_n\otimes \omega_A)((YY^*)^2)
&=&\sum\limits_{i, j=1}^d \nu_j(1-\nu_i)y_iy_j^*(y_iy_j^*)^*+\left(\sum\limits_{i=1}^d (1-\nu_i) y_iy_i^*\right)^2\label{eq20202020}\\
&=&\sum\limits_{i=1}^d (1-\nu_i)y_i\left(\sum\limits_{j=1}^d \nu_jy_j^*y_j\right)y_i^*+ \left(\sum\limits_{i=1}^d (1-\nu_i) y_iy_i^*\right)^2\nonumber\\
&\leq & \left(\sum\limits_{i=1}^d (1-\nu_i)y_iy_i^*\right)
\left(\left\|\sum\limits_{i=1}^d \nu_iy_i^*y_i\right\|+\left\|\sum\limits_{i=1}^d (1-\nu_i)y_iy_i^*\right\|\right)\nonumber\\
&=&\left(\left\|\sum\limits_{i=1}^d \nu_iy_i^*y_i\right\|+\left\|\sum\limits_{i=1}^d (1-\nu_i)y_iy_i^*\right\|\right)(\text{Id}_n\otimes \omega_A)(YY^*)\,,\nonumber
\end{eqnarray}
where the last equality is given by (\ref{eq20})\,.

As before, the crucial point is to show the following
\begin{lemma}\label{lem4242}
Let $x_1\,, \ldots \,, x_d\in
M_n(\mathbb{C})$\,. There exists $X\in M_n(\A)$ so that
$(\text{Id}_n\otimes E_A)(X)=\sum\limits_{i=1}^d x_i\otimes e_i$\,,
satisfying, moreover,
\begin{eqnarray}\label{eq7003}
\|X\|_{M_n(\A)}&\leq & \sqrt{2}\max\left\{\left\|\sum\limits_{i=1}^d
\nu_ix_i^*x_i\right\|^\frac12\,, \left\|\sum\limits_{i=1}^d
(1-\nu_i)x_ix_i^*\right\|^\frac12\right\}\,.
\end{eqnarray}
\end{lemma}

For this, we first prove the following

\begin{lemma}\label{lem1}
If $y_1\,, \ldots \,,y_d\in M_n(\mathbb{C})$ satisfy
\begin{eqnarray}\label{eq43}
\max\left\{\left\|\sum\limits_{i=1}^d
\nu_iy_i^*y_i\right\|^\frac12\,, \left\|\sum\limits_{i=1}^d
(1-\nu_i)y_iy_i^*\right\|^\frac12\right\}&=&1\,,
\end{eqnarray}
then there exists $Z\in M_n(\A)$ such that $\|Z\|_{M_n(\A)}\leq
\frac1{\sqrt{2}}$\,, and, moreover, when $z_1\,, \ldots \,,z_d$ are
defined by $(\text{Id}_n\otimes E_A)(Z)=\sum\limits_{i=1}^d
z_i\otimes e_i$\,, then
\begin{eqnarray*}
\max\left\{\left\|\sum\limits_{i=1}^d \nu_i
(y_i-z_i)^*(y_i-z_i)\right\|^\frac12\,, \left\|\sum\limits_{i=1}^d
(1-\nu_i)(y_i-z_i)(y_i-z_i)^*\right\|^\frac12\right\}&\leq
&\frac12\,.
\end{eqnarray*}
\end{lemma}

\begin{proof}
Set
\[ Y:=\sum\limits_{i=1}^d y_i\otimes a_i\in M_n(\A)\,. \]
Now let $C>0$ and define $F_C:\mathbb{R}\rightarrow \mathbb{R}$ by
formula (\ref{eq7070709822})\,. Use functional calculus to define
$Z\in M_n(\A)$ by (\ref{eq34345656560})\,. Then $\|Z\|_{M_n(\A)}\leq
C$ and, as shown in the proof of Lemma \ref{lem4567}\,, it follows
that
\[ (Y-Z)^*(Y-Z)\leq \frac1{16C^2}(Y^*Y)^2\,, \quad (Y-Z)(Y-Z)^*\leq \frac1{16C^2}(YY^*)^2\,. \]
By letting $z_i=(\text{Id}_n\otimes \phi_i^A)(Z)$\,, $1\leq i\leq d$\,, we then have
\[ (\text{Id}_n\otimes E_A)(Z)=\sum\limits_{i=1}^d z_i\otimes e_i\,,
\quad \text{respectively}\,, \quad (\text{Id}_n\otimes E_A)(Y)=
\sum\limits_{i=1}^d y_i\otimes E_A(a_i)=\sum\limits_{i=1}^d
y_i\otimes e_i\,,
\]
and we obtain the estimates
\begin{eqnarray*}
\sum\limits_{i=1}^d \nu_i(y_i-z_i)^*(y_i-z_i)&\leq &
(\text{Id}_n\otimes \omega_A)((Y-Z)^*(Y-Z))\\
&\leq & \frac1{16C^2}(\text{Id}_n\otimes \omega_A)((Y^*Y)^2)\\
&\leq & \frac1{16C^2}\left(\left\|\sum\limits_{i=1}^d \nu_iy_i^*y_i\right\|+\left\|\sum\limits_{i=1}^d (1-\nu_i)y_iy_i^*\right\|\right)(\text{Id}_n\otimes \omega_A)(Y^*Y)\\
&\leq & \frac{2}{16C^2}(\text{Id}_n\otimes
\omega_A)(Y^*Y)\\&=&\frac1{8C^2}\sum\limits_{i=1}^d \nu_iy_i^*y_i\,,
\end{eqnarray*}
respectively, $\sum\limits_{i=1}^d (1-\nu_i)(y_i-z_i)(y_i-z_i)^*\leq
\frac1{8C^2}\sum\limits_{i=1}^d (1-\nu_i)y_iy_i^*$\,. We deduce that
\begin{eqnarray*}\left\|\sum\limits_{i=1}^d \nu_i(y_i-z_i)^*(y_i-z_i)\right\|&\leq &\frac1{8C^2}\left\|\sum\limits_{i=1}^d \nu_iy_i^*y_i\right\|\,\,\leq \,\,\frac1{8C^2}\,,\end{eqnarray*}
respectively,
\begin{eqnarray*}
\left\|\sum\limits_{i=1}^d (1-\nu_i)(y_i-z_i)(y_i-z_i)^*\right\|&\leq & \frac1{8C^2}\left\|\sum\limits_{i=1}^d(1-\nu_i)y_iy_i^*\right\|\,\,\leq \,\,\frac1{8C^2}\,. \end{eqnarray*}
Hence
\begin{eqnarray*}
\max\left\{\left\|\sum\limits_{i=1}^d \nu_i
(y_i-z_i)^*(y_i-z_i)\right\|^\frac12\,, \left\|\sum\limits_{i=1}^d
(1-\nu_i)(y_i-z_i)(y_i-z_i)^*\right\|^\frac12\right\}&\leq
&\frac1{\sqrt{8}C}\,.
\end{eqnarray*}
Now take $C=\frac1{\sqrt{2}}$ to obtain the conclusion.
\end{proof}

We are now ready to prove Lemma \ref{lem4242}\,. Indeed, Lemma
\ref{lem1} shows that if $y:= \sum\limits_{i=1}^d y_i\otimes e_i\in
M_n(H)$ has norm $\|y\|_{M_n(H)}=1$\,, then there exists $Z\in
M_n(\A)$ such that $\|Z\|_{M_n(\A)}\leq \frac1{\sqrt{2}}$ and
$\|(\text{Id}_n\otimes E_A)(Z)-y\|_{M_n(H)}\leq  \frac12$\,. By
homogeneity we infer that for all $y\in M_n(H)$\,, there exists
$Z\in M_n(\A)$ satisfying the conditions $\|Z\|_{M_n(\A)}\leq
\frac1{\sqrt{2}}\|y\|_{M_n(H)}$ and $\|(\text{Id}_n\otimes
E_A)(Z)-y\|_{M_n(H)}\leq {\frac12}\|y\|_{M_n(H)}$\,. Applying now
Lemma \ref{lem25} with $C=\frac1{\sqrt{2}}$ and $\delta=\frac12$ we
deduce that for all $x\in M_n(H)$ there exists $X\in M_n(\A)$ so
that $(\text{Id}_n\otimes E_A)(X)=x$\,, satisfying, moreover,
\begin{eqnarray*}
\|X\|_{M_n(\A)}&\leq &
\frac{C}{1-\delta}\|x\|_{M_n(H)}=\sqrt{2}\|x\|_{M_n(H)}\,.
\end{eqnarray*}
The proof of Lemma \ref{lem4242} is complete.

By Lemmas \ref{lem4141} and \ref{lem4242} and Remark \ref{rem5}\,,
there exists a linear bijection $\widetilde{E}_A:
{{\A}/{\text{Ker}(E_A)}}\rightarrow H$ such that
\begin{eqnarray*}
\widetilde{E}_A (q_A(a_i))&=&e_i\,, \quad 1\leq i\leq d\,,
\end{eqnarray*}
making the following diagram commutative:
$$
\xymatrix{
 {M_n(\A)}\ar@{->}^{\text{Id}_n\otimes E_A}[rr]
 \ar@{->}_{\text{Id}_n\otimes q_A}[dr]
 & & {M_n(H)}\\
& {M_n({\A}/{\text{Ker}(E_A)})}
\ar@{->}[ur]_{\text{Id}_n\otimes \widetilde{E}_A}}
$$
Moreover, with respect to the natural operator space structure of
the quotient ${{\A}/{\text{Ker}(E_A)}}$ one has for all $x_1\,,
\ldots \,,x_d\in M_n(\mathbb{C})$\,,
\begin{eqnarray*}
\max\left\{\left\|\sum\limits_{i=1}^d
\nu_ix_i^*x_i\right\|^\frac12\,, \left\|\sum\limits_{i=1}^d
(1-\nu_i)x_ix_i^*\right\|^\frac12\right\}&\leq &
\left\|\sum\limits_{i=1}^d x_i\otimes
q_A(a_i)\right\|_{M_n({{\A}/{\text{Ker}(E_A)}})}\\
&\leq & \sqrt{2}\max\left\{\left\|\sum\limits_{i=1}^d
\nu_ix_i^*x_i\right\|^\frac12\,, \left\|\sum\limits_{i=1}^d
(1-\nu_i)x_ix_i^*\right\|^\frac12\right\}\,,
\end{eqnarray*}
i.e., the inequalities (\ref{eq10}) hold. This completes the proof
of Theorem \ref{th1} in the finite dimensional case.

We now consider the infinite dimensional case
($\text{dim}(H)=\infty$)\,. Let $V\subset H$ be a finite dimensional
subspace, and let $d=\text{dim}(V)$\,. Set
\[ A_V:={P_VA}_{\vert_{P_VH}}\in {\mathcal B}(P_VH)\,, \]
where $P_V$ is the projection of $H$ onto $V$\,. Then $0\leq A_V\leq
I$\,.

Let $\A_V$ be the CAR algebra on $V$\,, and denote by $\omega_A$
(respectively, $\omega_{A_V}$) the gauge-invariant quasi-free state
on ${\mathcal A}$ (respectively, $\A_V$) corresponding to the
operator $A$ (respectively, $A_V$)\,. Recall that $\A_V$ is the norm
closure of $\,{\text{Span}\{{a(e_{i_1})}^*\ldots
{a(e_{i_n})}^*{a(e_{j_1})}\ldots (e_{j_m})\,; \,1\leq i_1, \ldots
\,,i_n, j_1, \ldots \,,j_m, n,m\leq d\}}$\,. By equation (\ref{eq6})
it follows that ${\omega_A}_{\vert_{\A_V}}$ and $\omega_{A_V}$
coincide on all polynomials that generate $\A_V$\,. Since states are
norm continuous, we conclude that
 \begin{eqnarray}\label{eq4000}
{\omega_A}_{\vert_{\A_V}}&=& \omega_{A_V}\,.
\end{eqnarray}
The key point that will allow us to reduce the infinite dimensional
case to the finite dimensional one is the fact, which we will
justify in the following, that $E_A(b)\in V$\,, whenever $b\in
{\A}_V$\,. Indeed, given $b\in {\A}_V$\,, we will show that
\[ \langle E_A(b), f\rangle_H=0\,, \quad \forall f\in
V^\perp\,. \] By (\ref{eq8}), this is equivalent to showing that
\begin{eqnarray}\label{eq78788}
\omega_A(b{a(f)}^*+{a(f)}^*b)&=&0\,, \quad \forall f\in V^\perp\,.
\end{eqnarray}
By continuity, it suffices to consider elements $b\in {\A}_V$ of the
form $b={a(e_{i_1})}^*\ldots {a(e_{i_n})}^*{a(e_{j_1})}\ldots
(e_{j_m})$ where $1\leq i_1, \ldots \,,i_n, j_1, \ldots \,,j_m,
n,m\leq d$\,. Let $f\in V^\perp$\,. Since $f\perp e_i$\,, $1\leq
i\leq d$\,, we get by the CAR relations (\ref{eq4}) that
\begin{eqnarray*}
b{a(f)}^*&=&(-1)^{n+m}{a(f)}^*b\,.
\end{eqnarray*}
So if $n+m$ is odd, then $b{a(f)}^*+{a(f)}^*b=0$\,. If $n+m$ is
even, i.e., $n+m+1$ is odd, then by (\ref{eq4}) and (\ref{eq6})
(together with (\ref{eq4000})) it follows that
$\omega_A(b{a(f)}^*)=0=\omega_A({a(f)}^*b)$\,. Hence, in both cases
(\ref{eq78788}) follows, and our claim is proved. By uniqueness in
the construction of the maps $E_A$ and $E_{A_V}$\,, we conclude that
\begin{eqnarray}\label{eq5001}
{E_A}_{\vert_{\A_V}}&=& E_{A_V}\,.
\end{eqnarray}
Since $\A$ is the $C^*$-algebra generated by the operators
$a(\xi)$\,, $\xi\in H$\,, it is clear that
\begin{equation}\label{eq3031323334}
\A=\overline{\bigcup_V {\A_V}}\,,\quad (\text{norm closure})
\end{equation}
where the union is taken over all finite dimensional subspaces $V$
of $H$\,. Moreover, note that $\A_{V_1}\subseteq  \A_{V_2}$ when
$V_1\subseteq V_2$\,. We also claim that
\begin{equation}\label{eq3031323344}
\text{Ker}(E_A)=\overline{\bigcup_V {\text{Ker}(E_{A_V})}}\,,\quad
(\text{norm closure})
\end{equation}
where the right-hand side is also an increasing union because
${\text{Ker}(E_{A_V})}=\text{Ker}(E_A)\cap {\A_V}$\,, for all
$V\subset H$\,, finite dimensional subspace. To prove
(\ref{eq3031323344}), let $b\in \text{Ker}(E_A)$ and choose $b_n\in
\bigcup_V {\A_V}$\,, $n\geq 1$ such that $\|b_n-b\|\rightarrow 0$ as
$n\rightarrow \infty$\,. Further, set
\[ b'_n:=b_n-a(E_A(b_n))\,, \quad n\geq 1\,. \]
For $n\geq 1$\,, since $b_n\in {\A_{V_n}}$ for some finite
dimensional subspace $V_n$ of $H$\,, we have by (\ref{eq5001}) that
$E_A(b_n)=E_{A_{V_n}}(b_n)\in V_n$, and hence $b'_n\in \A_{V_n}$\,.
Moreover, by (\ref{eq9})\,, we get
\[ E_A(b'_n)=E_A(b_n)-E_A(b_n)=0\,. \]
Therefore, $b'_n\in \text{Ker}(E_A)\cap
\A_{V_n}=\text{Ker}(E_{A_{V_n}})$\,, which proves
(\ref{eq3031323344})\,. Now, since the union in formula
(\ref{eq3031323344}) is increasing, we also have for all $n\in
\mathbb{N}$\,,
\begin{equation}\label{eq3031323347}
M_n(\text{Ker}(E_A))=\overline{\bigcup_V
M_n(\text{Ker}(E_{A_V}))}\,.\quad (\text{norm closure})
\end{equation}
We are now ready to proceed with the proof of Theorem \ref{th1} in
the case $\text{dim}(H)=\infty$\,. We shall prove that for all
positive integers $n, r$ and all $x_i\in M_n(\mathbb{C})$ and
$b_i\in \A=\A(H)$\,, $1\leq i\leq r$\,, the inequalities
(\ref{eq13}) hold.

Indeed, by (\ref{eq3031323334}) and the fact that $\bigcup_V \A_V$
is an increasing union, it suffices to prove (\ref{eq13}) for
elements $b_i\in \A_{V_0}$\,, where $V_0$ is an arbitrary finite
dimensional subspace of $H$\,. Let now such $V_0$ be fixed. Since
Theorem \ref{th1} has been proved in the final dimensional case, we
have for each finite dimensional subspace $V$ with $V_0\subseteq
V\subset H$ that
\begin{equation}\label{eq30313233477}
\left\|\sum\limits_{i=1}^r x_i\otimes
E_{A_V}(b_i)\right\|_{M_n(V)}\leq \left\|\sum\limits_{i=1}^d
x_i\otimes q_A(b_i)\right\|_{M_n({\A_V}/{\text{Ker}(E_{A_V})})}\leq
\sqrt{2}\left\|\sum\limits_{i=1}^r x_i\otimes
E_{A_V}(b_i)\right\|_{M_n(V)}\,.
\end{equation}
By (\ref{eq5001})\,, $E_{A_V}(b_i)=E_A(b_i)$\,, $1\leq i\leq r$\,,
for all such $V$\,, since $b_i\in V_0\subseteq V$\,. Moreover, since
the norm in $M_n({\A}/{\text{Ker}(E_A)})$ is the quotient norm of
the quotient space ${M_n(\A)}/{M_n(\text{Ker}(E_A))}$\,, and
likewise for $M_n({\A_V}/{\text{Ker}(E_{A_V})})$\,, we get by
(\ref{eq3031323347}) that
\[ \lim\limits_V \left\|\sum\limits_{i=1}^r x_i\otimes
q_A(b_i)\right\|_{M_n({\A_V}/{\text{Ker}(E_{A_V})})}=\left\|\sum\limits_{i=1}^r
x_i\otimes b_i\right\|_{M_n({\A}/{\text{Ker}(E_A)})}\,, \] where the
limit is taken over the directed set of finite dimensional subspaces
$V$ with $V_0\subseteq V\subset H$\,, ordered by inclusion. Hence,
the inequalities (\ref{eq13}) follow from (\ref{eq30313233477}) and
the proof of Theorem \ref{th1} is complete.

\begin{cor}\label{ohembed}
Let $P$ be the hyperfinite type III$_1$ factor. For any subspace $H$
of $R\oplus C$\,, its dual $H^*$ embeds completely isomorphically
into the predual $P_*$ of $P$\,, with cb-isomorphism constant $\leq
\sqrt{2}$\,. In particular, the operator Hilbert space $OH$
cb-embeds into $P_*$ with cb-isomorphism constant $\leq \sqrt{2}$\,.
\end{cor}
\begin{proof}
Given a subspace $H$ of $R\oplus C$, let $A$ be the associated
operator ($0\leq A\leq I$) defined by (\ref{eq111}), ${\mathcal A}$
the CAR algebra over $H$, and $\omega_A$ the corresponding
gauge-invariant quasi-free state on $\A$\,. Denote
$\overline{\pi_A(\A)}^{\text{sot}}$ by $M$\,, where $\pi_A$ is the
unital $*$-homomorphism from the GNS representation associated to
$({\mathcal A}, \omega_A)$\,. By Theorem 5.1 in \cite{PS}, $M$ is a
hyperfinite factor. Then the von Neumann algebra tensor product
$M\bar{\otimes} P$ is (isomorphic to) the hyperfinite type III$_1$
factor $P$. Moreover, $M_*$ cb-embeds into $(M\bar{\otimes} P)_*$\,,
the embedding being given by the dual map of a normal conditional
expectation from $M\bar{\otimes} P$ onto $M$\,. Therefore, by
Theorem \ref{th1} it follows that the dual $H^*$ of $H$ embeds
completely isomorphically into $P_*$\,, with cb-isomorphism constant
$\leq \sqrt{2}$\,. Furthermore, note that $H^*$ is completely
isometric to a quotient of the dual space $(R{\oplus}_\infty C)^*$.
We infer that any quotient (and further, any sub-quotient, that is,
subspace of a quotient) of $(R{\oplus}_\infty C)^*$ cb-embeds into
$P_*$, with cb-isomorphism constant $\leq \sqrt{2}$\,. As shown by
Pisier (cf. \cite{Pi2}), the operator space $OH$ is a subspace of a
quotient of $R{\oplus}_\infty C$\,. Since $OH$ is self-dual as an
operator space (cf. \cite{Pi6})\,, $OH$ is also a sub-quotient of
$(R\oplus_\infty C)^*$\,. We conclude that $OH$ embeds completely
isomorphically into $P_*$, with cb-isomorphism constant $\leq
\sqrt{2}$\,. (See also Junge's results in Section 8 of \cite{Ju} on
the embedding of $OH$ into $P_*$\,.)
\end{proof}

\begin{rem}\label{remdualnorm}\rm
Let $H\subset R\oplus C$ be a subspace of dimension $d< \infty$\,,
and let $A$ be the associated operator defined by (\ref{eq111}),
respectively, let $\{e_i\}_{i=1}^d$\,, $\{\nu_i\}_{i=1}^d$ be
defined by (\ref{eq77777777789778}). Assume further that $0< \nu_i<
1$\,, for all $1\leq i\leq d$\,. Now define for any $x_1\,,
\ldots\,, x_d\in M_n(\mathbb{C})$\,,
\begin{equation}\label{eq3000}
|||\{x_i\}_{i=1}^d|||^*:=\inf\left\{\text{Tr}\left[\left(\sum\limits_{i=1}^d
\frac1{\nu_i} v_iv_i^*\right)^\frac12+ \left(\sum\limits_{i=1}^d
\frac1{(1-\nu_i)} z_i^*z_i\right)^\frac12\right]\,; \,x_i=v_i+z_i\in
M_n(\mathbb{C})\right\},
\end{equation}
where $\text{Tr}$ denotes, as before, the non-normalized trace on
$M_n(\mathbb{C})$\,.

Note that $|||\cdot|||^*$ is the dual norm of $\|\cdot\|_{M_n(H)}$.
From the proof of Theorem \ref{th1} it follows by duality that the
transpose $F_A:=E_A^*$ of the map $E_A:\A\rightarrow H$ becomes a
complete injection of $H^*$ into $\text{Span}\{\phi_1^A\,, \dots\,,
\phi_d^A\}={\A}^*$\,. More precisely, we obtain that
\begin{equation}\label{eq1000}
\frac1{\sqrt{2}}|||\{x_i\}_{i=1}^d|||^*\,\,\leq
\,\,\left\|\sum\limits_{i=1}^d x_i\otimes
\phi_i^A\right\|_{{M_n({\A})}^*}\leq |||\{x_i\}_{i=1}^d|||^*\,.
\end{equation}
\end{rem}

We now discuss estimates for best constants in the inequalities
(\ref{eq1000}) above.

\begin{theorem}\label{th789542}
Let $c_1\,, c_2$ denote the best constants in the inequalities
\begin{eqnarray}\label{eq545454221}
c_1|||\{x_i\}_{i=1}^d|||^*&\leq & \left\|\sum\limits_{i=1}^d
x_i\otimes \phi_i^A\right\|_{{M_n({\A})}^*}\leq
c_2|||\{x_i\}_{i=1}^d|||^*\,.
\end{eqnarray}
where $d, n$ are arbitrary positive integers, $H\subseteq R\oplus C$
is a Hilbert space of dimension $\text{dim}(H)=d$ with associated
operator $A$ given by (\ref{eq111})\,, $\A$ is the CAR-algebra over
$H$, $\phi_1^A\,, \ldots \,, \phi_d^A$ are defined by (\ref{eq88}),
and $x_1\,, \ldots \,, x_d\in M_n(\mathbb{C})$\,. Then
\begin{equation}\label{eq7868900003}
c_1=\frac1{\sqrt{2}}\,, \qquad c_2=1\,.
\end{equation}
\end{theorem}

\begin{proof}
By (\ref{eq1000}) we obtain immediately the following estimates
\begin{equation}\label{eq45345465}
\frac1{\sqrt{2}}\leq c_1\leq c_2\leq 1\,.
\end{equation}

Next we prove that $c_1=\frac1{\sqrt{2}}$\,. Take $n=1$\,, $d=1$\,,
in which case $H=\mathbb{C}$\,, $A={\frac12}I_H$ and
$\A=M_2(\mathbb{C})$\,, and let
$x_1=I_{M_1(\mathbb{C})}=I_{\mathbb{C}}$\,. Then
$\phi_1^A(b)=\text{Tr}(a_1^*b)$\,, $\forall b\in \A$\,, where
$a_1=\left(\begin{array}
[c]{cc}%
0 & 1\\
0 & 0
\end{array} \right)\in \A$\,. Since $|a_1|$ is a projection with
$\text{Tr}(|a_1|)=1$\,, we get
$\|\phi_1^A\|_{\A^*}=\|a_1^*\|_{L^1(\A\,,
\text{Tr})}=\|a_1\|_{L^1(\A\,,
\text{Tr})}=\|\,\,|a_1|\,\,\|_{L^1(\A\,, \text{Tr})}=1\,.$ It is
easily checked by the definition (\ref{eq3000}) that
$|||x_1|||^*=\sqrt{2}$\,, hence, $\frac1{\sqrt{2}}|||x_1|||^*=
1=\|\phi_1^A\|_{\A^*}=\|x_1\otimes \phi_1^A\|_{{\A}^*}$. It follows
that $c_1\leq \frac1{\sqrt{2}}$\,, which together with
(\ref{eq45345465}) imply that $c_1=\frac1{\sqrt{2}}$\,.

We now prove that $c_2=1$\,. For this, given $d\in \mathbb{N}$\,,
let $H=\text{Span}\{e_{1i}\oplus e_{i1}; 1\leq i\leq d\}\subseteq
R\oplus C$\,. It follows easily by (\ref{eq111}) that the associated
operator is $A={\frac12}I_H$\,. Let $\{e_i\}_{i\geq 1}$ be an
orthonormal basis of $H$ with respect to which the matrix $A$ is
diagonal. As before, let \[ a_i:=a(e_i)\,, \quad 1\leq i\leq d \] be
the generators of the CAR algebra $\A=\A(H)$ built on $H$\,. We
consider the special representation of $\A$ constructed in the proof
of Lemma \ref{lem677767776} and use the identification
\begin{equation}\label{eq666666666444444409} \A\cong
\psi(\A)=\otimes_{i=1}^d M_2(\mathbb{C})\,, \end{equation} where
$\psi$ is the $*$-isomorphism obtained therein. Via this
identification, we may assume that the generators $a_i$\,, $1\leq
i\leq d$, of $\A$ are given by (\ref{eq300004509}). Note also that
the eigenvalues of $A$ are $\nu_i=\frac12$\,, $\forall 1\leq i\leq
d$\,, so the corresponding quasi-free state $\omega_A$ on $\A$ is
tracial. For simplicity of notation, let $\omega_A$ be denoted by
$\tau$\,.

For all $1\leq i\leq d$\,, set $x_i:=e_{i1}\in M_d(\mathbb{C})$\,.
In what follows, $\text{Tr}$ denotes the non-normalized trace on
$M_d(\mathbb{C})$\,. For $1\leq i\leq d$, we have
$\phi_i^A(b)=\tau(a_i^*b+ba_i^*)=2\tau(a_i^*b)$\,, $\forall b\in
\A$\,.

Let $h_i:=2a_i^*$\,, $1\leq i\leq d$\,. Then
\begin{eqnarray*}\left\|\sum\limits_{i=1}^d x_i\otimes
\phi_i^A\right\|_{(M_d({\A}))^*}&=&\left\|\sum\limits_{i=1}^d
x_i\otimes h_i\right\|_{L^1(M_d(\mathbb{C})\otimes \A,
\text{Tr}\otimes \tau)}\\&=& (\text{Tr}\otimes
\tau)\left(\left[\left(\sum\limits_{i=1}^d
x_i^*h_i\right)^*\left(\sum\limits_{i=1}^d
x_i^*h_i\right)\right]^\frac12\right)\\
&=&\tau\left(\left(\sum\limits_{i=1}^d
h_i^*h_i\right)^\frac12\right)=2\,\tau\left(\left(\sum\limits_{i=1}^d
a_ia_i^*\right)^\frac12\right)\,.
\end{eqnarray*}
Note that by (\ref{eq3000}) it follows immediately that
$|||\{x_i\}_{i=1}^d|||^*\leq
\text{Tr}\left(\left(2\sum\limits_{i=1}^d
x_i^*x_i\right)^\frac12\right)\leq \sqrt{2d}$\,. Therefore, if we
show that
\begin{equation}\label{eq888888888765}
\lim\limits_{d\rightarrow \infty}
\sqrt{\frac{2}{d}}\,\,\tau\left(\left(\sum\limits_{i=1}^d
a_ia_i^*\right)^\frac12\right)=1\,,
\end{equation} it then follows by (\ref{eq545454221}) that $c_2\geq 1$,
which implies that $c_2=1$\,.

We now prove (\ref{eq888888888765}). For this, we first show that
$a_1a_1^*\,, \ldots , a_da_d^*$ are independent, self-adjoint random
variables with distribution
\begin{equation}\label{eq55555555555333335}\mu_{a_ia_i^*}=\frac12(\delta_{\{0\}}+\delta_{\{1\}})\,, \quad
1\leq i\leq d\,. \end{equation} Using the notation set forth in the
proof of Lemma \ref{lem677767776} and (\ref{eq300004509}), a simple
computation shows that
\begin{equation}\label{eq55555555555333337}
a_1a_1^*=\left(\begin{array}
[c]{cc}%
1 & 0\\
0 & 0
\end{array}
\right)\otimes (\otimes_{j=2}^d I_2)\,, \quad
a_ia_i^*:=\left(\otimes_{j=1}^{i-1} I_2\right)\otimes
\left(\begin{array}
[c]{cc}%
1 & 0\\
0 & 0
\end{array}
\right)\otimes \left(\otimes_{j=i+1}^d I_2\right)\,, \quad 2\leq
i\leq d\,.
\end{equation}
In particular, $a_ia_i^*$ is a projection, for all $1\leq i\leq
d$\,. So $a_ia_i^*$ has spectrum $\sigma(a_ia_i^*)=\{0, 1\}$\,, and
since $\tau(a_ia_i^*)=\frac12$\,, formula
(\ref{eq55555555555333335}) follows.

By (\ref{eq55555555555333337}), $a_ia_i^*$ and $a_ja_j^*$ do
commute, for all $1\leq i, j\leq d$\,. Thus, in order to prove the
independence of $a_1a_1^*\,, \ldots , a_da_d^*$ (both in the
classical sense and in the sense of Voiculescu (cf. \cite{Vo})), it
remains to show that
\begin{equation*} \tau\left(
(a_1a_1^*)^{m_1}\ldots (a_da_d^*)^{m_d}\right)=\prod\limits_{i=1}^d
\tau((a_ia_i^*)^{m_i})\,, \quad m_1\,,\ldots, m_d\in \mathbb{N}\,.
\end{equation*}
This follows immediately from the special form
(\ref{eq55555555555333337}) of the elements $a_ia_i^*$\,, $1\leq
i\leq d$\,, and the fact that by the identification
(\ref{eq666666666444444409})\,, $\tau$ can be viewed as the tensor
product trace on $\otimes_{i=1}^d M_2(\mathbb{C})$\,.

Now recall that $d$ was arbitrarily chosen. By applying the law of
large numbers we deduce that the sequence
$\left\{{\frac1d}{\sum\limits_{i=1}^d a_ia_i^*}\right\}_{d\geq 1}$
converges in probability to ${\frac12}I_{M_n(\mathbb{C})}$\,, as
$d\rightarrow \infty$\,. This implies that
\begin{equation}\label{eq66666666644444440789}
\lim\limits_{d\rightarrow \infty}\sqrt{{\frac1d}{\sum\limits_{i=1}^d
a_ia_i^*}}={\frac1{\sqrt{2}}}I_{M_n(\mathbb{C})} \quad\text{in
probability}\,.
\end{equation}
Since, moreover, $0\leq {\frac1d}{\sum\limits_{i=1}^d a_ia_i^*}\leq
1$\,, for all $d\geq 1$\,, it follows that the convergence
(\ref{eq66666666644444440789}) holds also in the 2-norm. Hence
$\lim\limits_{d\rightarrow \infty} \tau\left(\left(\frac1d
\sum\limits_{i=1}^d
a_ia_i^*\right)^\frac12\right)=\frac1{\sqrt{2}}$\,, which gives
(\ref{eq888888888765}), and the proof is complete.
\end{proof}

\section*{Acknowledgements} We would like to thank Gilles Pisier for
his careful reading of a previous version of the manuscript, and his
very useful comments and suggestions, which improved the paper. In
particular, the argument leading to the sharp estimate $c_2=1$ in
Theorem \ref{th8998} and Theorem \ref{th789542} was suggested to us
by him.

\thanks{}

\end{document}